\documentclass[twoside,11pt]{article}
\title{} \author{} \date{}
\usepackage{latexsym,amssymb,times}
\input amssym.def
\newtheorem{te}{Theorem}[section]

\newtheorem{fac}[te]{Fact}
\newtheorem{cla}[te]{Claim}

\newtheorem{rem}[te]{Remark}
\newtheorem{ex}[te]{Example}

\def\dok{\noindent{\bf Proof. }}
\def\kdok{\hfill $\Box$ \par \vspace*{2mm} }
\def\a{\alpha}
\def\b{\beta}

\def\f{\varphi}
\def\p{\psi}
\def\o{\omega}
\def\k{\kappa}

\def\r{\rho}
\def\s{\sigma}

\def\t{\tau}

\def\E{{\mathbb E}}
\def\F{{\mathbb F}}
\def\S{{\mathbb S}}
\def\P{{\mathbb P}}
\def\Q{{\mathbb Q}}

\def\N{{\mathbb N}}
\def\X{{\mathbb X}}
\def\Y{{\mathbb Y}}
\def\Z{{\mathbb Z}}

\def\C{{\mathbb C}}
\def\BG{{\mathbb G}}

\def\BH{{\mathbb H}}
\def\BK{{\mathbb K}}

\def\CP{{\mathcal P}}

\def\CC{{\mathcal C}}

\def\CT{{\mathcal T}}

\def\K{{\mathcal K}}
\def\L{{\mathcal L}}

\def\CN{{\mathcal N}}
\def\CX{{\mathcal X}}

\def\CF{{\mathcal F}}
\def\CG{{\mathcal G}}

\def\la{\langle}
\def\ra{\rangle}

\def\deg{\mathop{\mathrm{deg}}\nolimits}

\def\id{\mathop{\mathrm{id}}\nolimits}

\def\Conv{\mathop{\rm Conv}\nolimits}

\def\Aut{\mathop{\rm Aut}\nolimits}

\def\Sym{\mathop{\rm Sym}\nolimits}

\def\Sent{\mathop{\rm Sent}\nolimits}
\def\Form{\mathop{\rm Form}\nolimits}

\def\Rado{\mathop{\mathrm{Rado}}\nolimits}
\def\ar{\mathop{\rm ar}\nolimits}
\def\Mod{\mathop{\rm Mod}\nolimits}

\def\Int{\mathop{\rm Int}\nolimits}
\def\Hom{\mathop{\rm Hom}\nolimits}

\def\At{\mathop{\rm At}\nolimits}

\def\Cond{\mathop{\rm Cond}\nolimits}

\def\Rev{\mathop{\rm Rev}\nolimits}
\def\sRev{\mathop{\rm sRev}\nolimits}
\def\wRev{\mathop{\rm wRev}\nolimits}

\def\Var{\mathop{\rm Var}\nolimits}
\def\Ord{\mathop{\rm Ord}\nolimits}
\def\Min{\mathop{\rm Min}\nolimits}
\def\Max{\mathop{\rm Max}\nolimits}
\def\Th{\mathop{\rm Th}\nolimits}

\begin{document}
\thispagestyle{plain}
\begin{center}
           {\large \bf \uppercase{Reversibility of extreme relational structures}}
\end{center}
\begin{center}
{\small\bf  Milo\v s S.\ Kurili\'c and Nenad Mora\v ca}\\
         Department of Mathematics and Informatics, University of Novi Sad, \\
         Trg Dositeja Obradovi\'ca 4, 21000 Novi Sad, Serbia.\\
                                     e-mail: milos@dmi.uns.ac.rs, nenad.moraca@dmi.uns.ac.rs
\end{center}
\begin{abstract}
\noindent
A relational structure
$\X$ is called reversible iff each bijective homomorphism from $\X$ onto $\X$ is an isomorphism, and linear orders are
prototypical examples of such structures. One way to detect new reversible structures of a given relational language $L$ is to notice that
the maximal or minimal elements of isomorphism-invariant sets of interpretations of the language $L$ on a fixed domain $X$ determine reversible structures.
We isolate certain syntactical conditions providing that a consistent $L_{\infty \o}$-theory defines a class of interpretations having extreme elements on a fixed domain and  detect several classes of reversible structures. In particular, we characterize the reversible countable ultrahomogeneous graphs.  
\end{abstract}
\section{Introduction}\label{S1}
Generally speaking, a structure is reversible iff each bijective endomorphism of that structure is an automorphism. Several prominent structures have this property;
for example, each compact Hausdorff space $\CX$ is reversible (because each continuous bijection $f:\CX \rightarrow \CX$ is a closed mapping and,
hence, a homeomorphism) and, similarly, each linear order $\X$ is a reversible relational structure
(since an increasing bijection $f:\X\rightarrow\X$ must be an isomorphism).

The reversible structures mentioned above are extreme: compact Hausdorff topologies are, on one hand, maximal compact and, on the other hand, minimal
Hausdorff topologies, and linear orders are maximal partial orders.
In this paper, searching for reversible structures, we investigate this phenomenon in the class of relational structures. So throughout the paper we assume that
$L=\la R _i:i\in I\ra$ is a relational language,  where $\ar (R_i)=n _i\in \N$, for $i\in I$, that $X$ is a non-empty set and
$\Int _L (X)=\prod _{i\in I}P(X^{n_i})$ the set of all interpretations of the language $L$, over the domain $X$. An interpretation
$\r= \la \r _i :i\in I\ra\in \Int _L (X)$ will be called reversible iff $\la X, \r\ra$ is a reversible structure.

First in Section \ref{S3} we easily establish the reversibility of minimal and maximal elements of the poset $\la \CC ,\subset \ra$, where
$\CC \subset \Int _L (X)$ is an isomorphism-invariant set, and, in particular, if $\CC $ is of the form
$\Int^\CT _L (X)=\{ \r \in \Int _L (X): \la X, \r\ra\models \CT\}$, for some set $\CT$ of sentences of the infinitary language $L_{\infty \o}$.
Of course, there are sets of the form $\Int^\CT _L (X)$ having neither minimal nor maximal elements, and, hence,
in Section \ref{S4} we isolate a class of formulas $\CF$
such that the set of maximal elements of the poset $\Int^\CT _L (X)$ is co-dense, whenever $\CT \subset \CF$,
and prove a dual statement concerning minimal elements.  We note that it is not our goal to find a syntactical characterization of
the largest class $\CF$ with the property mentioned above, because, for example,
for a countable language $L$, each isomorphism-invariant set $\CC \subset \Int _L (\o )$ is of the form $\Int^{\{ \f \} } _L (\o)$, where
$\f$ is the disjunction of the Scott sentences of the structures belonging
to $\CC$ and, trivially, the set $\Int^{\{ \f \lor \f _m \} } _L (\o)$, where
$\f _m:=\bigwedge _{i\in I} \forall \bar v \; R_i(\bar v)$ has a largest element, $\la X^{n_i} :i\in I\ra$.
Our goal is to find a reasonable class of sentences providing relevant examples of reversible structures.

Sections \ref{S6} and \ref{S7} contain some applications of the results mentioned above. In particular it is shown that the concept of ``forbidden finite substructures" provides a large class of extreme (and, hence, reversible) structures.
Clearly, one thing is to prove that extreme interpretations exist and the other is to find (or characterize) them.
Some results on this topic are given in examples.
\section{Preliminaries}\label{S2}
\paragraph{The algebra of interpretations}
Abusing notation, for $\r ,\s \in \Int _L (X)$ we will write
$\r\subset \s $ iff $\r_i \subset \s _i$, for all $i\in I$.
Clearly
$\la \Int _L (X), \subset \ra$ is a Boolean lattice and,
abusing notation again, the operations  in the corresponding Boolean algebra will be denoted in the following way:
if $\r ^j\in \Int _L (X)$, for $j\in J$, then
$\bigcap _{j\in J}\r ^j :=\la \bigcap _{j\in J}\r ^j_i:i\in I\ra$,
$\bigcup _{j\in J}\r ^j :=\la \bigcup _{j\in J}\r ^j_i:i\in I\ra$,
$\r ^c :=\la X ^{n_i}\setminus \r _i :i\in I\ra $,
$0 := \la \emptyset :i\in I\ra$ and
$1 := \la X^{n_i} :i\in I\ra$.
\paragraph{Direct and inverse images of interpretations}
If $X$ and $Y$ are non-empty sets and $n\geq 2$, the {\it n-th power} of a mapping $f: X\rightarrow Y$ is the mapping $f^n :X^n \rightarrow Y^n$  defined by
$f^n(\la  x_1, \dots , x_n\ra)= \la f(x_1), \dots , f(x_n) \ra $,
for each $\la  x_1, \dots , x_n\ra\in X^n$.
Clearly, $f$ is an injection (surjection) iff $f^n$ is an injection (surjection).

For $L$-interpretations $\r= \la \r _i :i\in I\ra\in \Int _L (X)$ and $\s= \la \s _i :i\in I\ra\in \Int _L (Y)$
the interpretations $f[\r ]\in \Int _L (Y)$ and $f^{-1}[\s ]\in \Int _L (X)$ are defined by
\begin{equation}\label{EQA000}
f[\r ]:= \la f^{n_i}[\r _i] :i\in I\ra \; \mbox{ and }\; f^{-1}[\s ]:= \la (f^{n_i})^{-1}[\s _i] :i\in I\ra,
\end{equation}
and these operators have all properties of direct and inverse images:
$\r \subset f^{-1}[f[\r ]]$,  $f[f^{-1}[\s]]= \s \cap f[1]$, $f[\bigcap _{j\in J}\r ^j]\subset \bigcap _{j\in J}f[\r ^j]$, $f^{-1}[\r ^c]=(f^{-1}[\r ])^c$, etc.
\paragraph{Morphisms}
Bijective homomorphisms will be called {\it condensations}.
Working with elements of $\Int _L (X)$, instead of $\Hom (\la X,\r \ra , \la X, \s \ra )$ we will write $\Hom (\r ,\s )$.
Also, instead of $\la X,\r \ra \cong \la X, \s \ra$ we will shortly write $\r \cong \s$ and regard
$\cong$ as an equivalence relation on the set $\Int _L (X)$.
Let $[\r ]_{\cong }:=\{ \s \in \Int _L (X): \s \cong \r \}$.
\begin{fac}\label{TA001}
For each $\r =\la \r _i :i\in I\ra , \s= \la \s _i :i\in I\ra \in \Int _L (X)$ we have:

(a) $\Hom (\r ,\s ) =\{ f\in {}^{X}X :f[\r ]\subset \s\}$;

(b) $[\r ]_{\cong }=\{ f[\r ]: f\in \Sym (X )\}$.
\end{fac}
\paragraph{The condensation order and reversibility}
If $\P =\la P , \leq \ra$ is a partial order,
a subset $C$ of $P$ is called
{\it convex} iff $p\leq q \le r$ and $p,r\in C$ implies $q\in C$. A set $A\subset P$ is called
an {\it antichain} iff different elements of $A$ are incomparable. Clearly, each antichain is convex and
$\Conv _\P (A) =\{ p\in P: \exists a',a'' \in A \;\; a' \leq p\leq a''\}$ is the minimal convex set containing the set $A\subset P$
(the {\it convex closure} of $A$).

Here we recall some facts from \cite{KMvar,KMsim,KM3}. Let $\preccurlyeq _c$ be the pre-order on the set $\Int _L (X)$
defined by: $\r \preccurlyeq _c \s $ iff there is a condensation $f: \la X , \r \ra \rightarrow \la X , \s \ra$.
The corresponding antisymmetric quotient, the poset $\la \Int _L (X) /\!\!\sim _c, \leq _c\ra$, where
$\r \sim _c \s \Leftrightarrow \r \preccurlyeq _c \s \land \s \preccurlyeq _c \r$ and
$[\r ]_{\sim _c}\leq  _c[\s ]_{\sim _c} \Leftrightarrow \r \preccurlyeq _c \s$, for $\r ,\s \in \Int _L (X)$,
is called the {\it condensation order}. Defining $[\r ]_{\sim _c }:=\{ \s \in \Int _L (X): \s \sim _c \r \}$, for each $\r \in \Int _L (X)$ we have
\begin{equation}\label{EQA027}
[\r ]_{\cong } \subset [\r ] _{\sim _c}=\Conv _{\la \Int _L (X), \subset \ra} ([\r ]_{\cong }).
\end{equation}
\begin{fac}\label{T4312}
For each interpretation $\r \in \Int_L(X)$ the following is equivalent:

(a) $\r$ is reversible, that is $\Cond (\r )=\Aut (\r )$,

(b) $[\r ]_{\cong}$ is an  antichain in the Boolean lattice $\la \Int_L(X),\subset \ra$,

(c) there is no $\s \in [\r ]_{\cong }$ such that $\r \varsubsetneq \s$,

(d) there is no $\s \in [\r ]_{\cong }$ such that $\s \varsubsetneq \r$,

(e) $\r ^c$ is reversible.
\end{fac}
An interpretation $\r \in \Int _L (X)$ will be called
{\it strongly reversible} iff $[\r ]_{\cong }=\{ \r \}$;
{\it weakly reversible} iff $[\r ]_{\cong }$ is a convex set in the Boolean lattice $\la \Int _L (X),\subset \ra$.
Clearly we have $\sRev_L(X)\subset \Rev_L(X)\subset \wRev_L(X)$, where
$\Rev_L(X)$ (resp.\ $\sRev_L(X)$, $\wRev_L(X)$) denotes the set of all
reversible (resp.\ strongly reversible, weakly reversible) interpretations $\r \in \Int _L (X)$.

It is easy to see that both reversibility and its two variations are $\sim _c$-invariants and, hence, $\cong$-invariants.
(A property $\CP$
is called {\it $\sim$-invariant} iff for each $\r,\s\in \Int _L (X)$
we have: if $\r$ has $\CP$ and $\s\sim \r$, then $\s$ has $\CP$).
In addition, weakly reversible interpretations have the Cantor-Schr\"{o}der-Bernstein
property for condensations (if $\r$ is  weakly reversible and  there are condensations $f:\la X, \s\ra \rightarrow \la X, \r\ra$ and
$g:\la X, \r\ra \rightarrow \la X, \s\ra$, then $\s\cong \r$).

Concerning strong reversibility we have: an interpretation $\r \in \Int _L (X)$ is strongly reversible iff for each $i\in I$,
the relation $\r _i $ is a subset of $X^{n_i}$ which is definable by a first-order formula of the empty language without parameters.
\begin{ex}\rm\label{EXA012}
If $L_b =\la R \ra$ is the binary language (i.e.\ $\ar(R)=2$) and $X\neq\emptyset$,
then the only strongly reversible elements of $\Int _{L_b}(X)$
are: $\emptyset$ (the empty relation), $\Delta _X$ (the diagonal), $X^2\setminus \Delta _X$ (the complete graph) and $X^2$ (the full relation).
\end{ex}
\paragraph{Partial orders}
If $\P =\la P , \leq \ra$ is a partial order,
by $\Min \P $ (resp.\ $\Max \P $) we denote the set of minimal (resp.\ maximal) elements of $\P$.
A set $D\subset P$ is called {\it dense} (resp.\  {\it co-dense}) in $\P$ iff for each $p\in P$ there is $q\in D$ such that $q\leq p$ (resp.\ $p\leq q$).
A set of $L$-interpretations $\CC \subset \Int _L (X)$ will be called {\it union-complete} (resp.\ {\it intersection-complete}) iff $\bigcup \L \in \CC$
(resp.\ $\bigcap \L \in \CC$) for each chain $\L \subset \CC$.
The partial order $\la \CC ,\subset \ra$ will be shortly denoted by $\CC$, when it is convenient.
\begin{fac}\label{TA028}
If $\CC \subset \Int _L (X)$ is a union-complete (resp.\ intersection-complete) set, then
$\Max \CC $  (resp.\ $\Min  \CC $) is a co-dense (resp.\  dense) subset of $\CC $.
\end{fac}
\dok
If $\CC $ is union-complete and $\r\in \CC$, then, by the Hausdorff maximal principle, there is a maximal chain $\L$ in $\CC$
such that $\r \in \L$. By the union-completeness of $\CC$, we have $\bigcup\L \in \CC$, by the maximality of $\L$ we have
$\bigcup\L \in \Max \la \CC ,\subset \ra  $ and, since $\r\in \L$, we have $\r\subset \bigcup\L$. The proof for intersection-complete sets is dual.
\hfill $\Box$
\paragraph{Infinitary languages}
Let $L =\la R_i :i\in I\ra$ be a relational language, $\k$ an infinite cardinal and $\Var =\{ v_\a : \a \in \k\}$ a set of variables.
By $\At _{L }$ we denote the corresponding set of {\it atomic formulas}, that is,
$$
\At _{L }=\{ v_\a =v_\b : \a ,\b \in \k\} \cup \{ R_i (v_{\a _1}, \dots , v_{\a _{n_i}}): i\in I \land  \la \a _1, \dots , \a _{n_i} \ra \in \k ^{n_i} \}.
$$
The class of {\it $L_{\infty  \o }$-formulas} is the class $\Form _{L_{\infty  \o }}= \bigcup _{\xi \in \Ord}\Form _\xi$, where
\begin{eqnarray*}
\Form _0        &\! =  \!  & \At _{L },\\
\Form _{\xi +1} &\! =  \!  & \Form _\xi  \; \cup \;  \{ \neg \f :  \f\in \Form _\xi \}\\
                &\! \cup \!& \{ \forall v_\a \; \f : \a\in \k \land  \f\in \Form _\xi  \} \; \cup  \; \{ \exists v_\a \; \f : \a\in \k \land \f\in \Form _\xi  \} \\
                &\! \cup \!& \textstyle \{ \bigwedge \Phi :  \Phi \subset \Form _\xi  \} \;  \cup  \; \{ \bigvee \Phi :  \Phi \subset \Form _\xi  \} , \\
\Form _{\gamma} &\! =   \! & \textstyle  \bigcup _{\xi <\gamma }\Form _\xi , \mbox{ for a limit ordinal } \gamma .
\end{eqnarray*}
Let $\X =\la X, \la R_i^{\X }:i\in I \ra$ be an $L$-structure and $\vec x =\la x_\a :\a \in \k \ra\in {}^{\k }X$ a valuation.
If $\b\in \k$ and $x'\in X$, by $\vec x _{\la \b , x'\ra}$ we denote the valuation $\vec y \in {}^{\k }X$ defined by:
$y _\a =x_\a$, for all $\a\neq \b$; and $y_\b =x'$.
The satisfiability relation for $L_{\infty  \o }$-formulas
is defined
in a standard way, for example,
$\X \models (R_i (v_{\a _1}, \dots , v_{\a _{n_i}}))[\vec x ]$ iff  $\la x_{\a _1}, \dots , x_{\a _{n_i}} \ra \in R_i^{\X }$;
$\;\X \models (\neg \f )[\vec x ]$ iff $\X \models \f [\vec x ]$ is not true;
$\;\X \models (\forall v_\a \; \f )[\vec x ]$ iff $\X \models \f [\vec x _{\la \a , x \ra} ]$, for each $x\in X$;
$\X \models (\bigvee \Phi )[\vec x ]$ iff $\X \models \f [\vec x ]$, for some $\f \in \Phi $.

$L_{\infty  \o }$-formulas $\f$ and $\p$ are called {\it logically equivalent}, in notation $\f \leftrightarrow \p$ iff they
are equivalent in all $L$-structures, that is, iff for each  $L$-structure $\X$ we have:
\begin{equation}\label{EQA036}
\forall \vec x \in {}^{\k } X \;\;\Big(\X \models \f [\vec x ] \Leftrightarrow \X \models \p [\vec x]\Big).
\end{equation}
If $\X$ and $\Y$ are $L$-structures, a mapping $f:X\rightarrow Y$ {\it preserves an $L_{\infty  \o }$-formula}
$\f$ iff
\begin{equation}\label{EQA021}
\forall \vec x \in {}^{\k } X \;\;\Big(\X \models \f [\vec x ] \Rightarrow \Y \models \f [f\vec x]\Big),
\end{equation}
where $f\vec x= \la f(x_\a ):\a \in \k \ra$. We say that the formula $\f$ is {\it absolute under} $f$ iff
in (\ref{EQA021}) we have ``$\Leftrightarrow$" instead of ``$\Rightarrow$".
\section{Reversibility of maximal and minimal structures}\label{S3}
A set $\CC\subset \Int _L (X)$
is said to be {\it isomorphism-invariant} or, shortly, {\it $\cong$-invariant } iff
\begin{equation}\label{EQA020}
\forall \r \in \CC \;\; [\r ]_{\cong }\subset \CC .
\end{equation}
By  $\CC  ^c $ we will denote the set $\{ \r ^c : \r \in \CC \}$. (Clearly, $\CC  ^c \neq \Int _L (X)\setminus \CC$.)
\begin{te}\label{TA018}
If $\CC\subset \Int _L (X)$ is an $\cong$-invariant set and
$\t\in \Max \CC$ (respectively, $\t\in \Min \CC$), then

(a) $\t$ is a reversible interpretation;

(b) $[\t ]_{\cong }=[\t ]_{\sim _c }$ is an antichain in $\CC$ and $[\t ]_{\cong }\subset \Max \CC $  (resp.\ $[\t ]_{\cong }\subset \Min  \CC  $);

(c) The set $\CC  ^c$ is $\cong$-invariant and $\t \in \Max \CC $ iff $\t ^c \in \Min (\CC  ^c )$, for any $\t$.
\end{te}
\dok
(a) Suppose that $\t$ is not reversible. Then, by Fact \ref{T4312}, there is $\s \in [\t ]_{\cong }$ such that $\t \varsubsetneq \s$
and, by (\ref{EQA020}), $\s \in \CC$, which is impossible, by the maximality of $\t$.

(b) By (a) and Fact \ref{T4312}, $[\t ]_{\cong }$ is an antichain and, by (\ref{EQA027}),  $[\t ]_{\cong }=[\t ]_{\sim _c }$.
Suppose that there are $\t _1\in [\t ]_{\cong }$ and $\r \in \CC$ such that $\t _1 \varsubsetneq \r$. Then, by Fact \ref{TA001}(b) there is
$f\in \Sym (X)$ such that $f[\t _1]=\t$, which, together with (\ref{EQA020}) implies $\t \varsubsetneq f[ \r ]\in \CC$. But,
by the maximality of $\t$, this is impossible. For $\t\in \Min \CC$ the proof is dual.

(c) For $\r\in \CC$ we show that $[\r ^c]_{\cong }\subset \CC  ^c$. So, if $\s \in [\r ^c]_{\cong }$, then, by Fact \ref{TA001}(b), there is
a bijection $f:X\rightarrow X$ such that $\s = f[\r ^c ]=f[\r ]^c$ and, since $f[\r ]\in [\r ]_{\cong }\subset \CC$, we have
$\s \in \CC  ^c$. Let $\t \in \Int _L (X)$. Then $\t ^c \in \Min (\CC  ^c )$ iff
$\t ^c =\r ^c$, for some $\r \in \CC$, and for each $\r \in \CC$ we have $\r ^c \subset \t ^c \Rightarrow \r ^c = \t ^c$.
In other words, $\t =\r $, for some $\r \in \CC$  and  for each $\r \in \CC$ we have $\t \subset \r \Rightarrow \r =\t$, which means that
$\t \in \Max \CC $.
\kdok
We note that, in fact, an interpretation $\t \in \Int _L (X)$ is reversible iff   $\t\in \Max \CC $
for some $\cong$-invariant set $\CC\subset \Int _L (X)$.
Namely, the remaining implication is trivial: if $\t$ is reversible, then by Fact \ref{T4312} we have
$\t\in [\t ]_{\cong}= \Max  [\t ]_{\cong} $.

Now we consider the sets of interpretations satisfying $L_{\infty  \o }$-sentences. In order to make a correspondence between
interpretations and their complements to each $L_{\infty  \o }$-formula $\f$ we adjoin the formula $\f  ^c$ defined in the following way:

 $(v_\a =v_\b ) ^c := v_\a =v_\b $ and $(R_i (v_{\a _1}, \dots , v_{\a _{n_i}})) ^c := \neg  R_i (v_{\a _1}, \dots , v_{\a _{n_i}}) $;

\noindent
If $\xi\in \Ord$ and $\f  ^c$ is defined for a formula $\f\in \Form _\xi$, then

$(\neg \f  ) ^c := \neg \f  ^c   $, $(\forall v_\a \; \f ) ^c := \forall v_\a \; \f  ^c  $ and
$(\exists v_\a \; \f ) ^c := \exists v_\a \; \f  ^c  $;

\noindent
If $\Phi \subset \Form _\xi$ and $\f  ^c$ is defined for each formula $\f\in \Phi$, then

$(\bigwedge \Phi  ) ^c :=\bigwedge \Phi  ^c  $ and $(\bigvee \Phi  ) ^c := \bigvee \Phi  ^c $,

\noindent
where, for a set $\Phi$ of $L_{\infty  \o }$-formulas $\Phi  ^c $ denotes the set $\{ \f  ^c : \f \in \Phi \}$.

\begin{te}\label{TA027}
If $X$ is a non-empty set and  $\CT$ a set of $L_{\infty  \o }$-sentences, then

(a) The set $\Int ^\CT _L (X)$ is $\cong$-invariant;

(b) $(\Int ^\CT _L (X)) ^c = \Int ^{\CT  ^c}_L (X) $ and this set is $\cong$-invariant;

(c) Maximal and minimal elements of $\Int ^\CT _L (X)$ are reversible interpretations;

(d) $\t\in \Max \Int ^\CT _L (X)$ iff
$\t ^c\in\Min \Int ^{ \CT  ^c  }_L (X)$, for  $\t \in \Int _L (X)$.
\end{te}
\dok
(a) If $\r \in \Int ^\CT _L (X)$ and $\s\in [\r ]_{\cong }$, then there exists an isomorphism $f:\la X, \r \ra \rightarrow \la X, \s \ra$.
Since for each sentence $\f \in \CT$ we have $\la X, \r \ra \models \f$, by a standard fact that
each $L_{\infty  \o }$-formula is absolute under each isomorphism
 we have $\la X, \s \ra \models \f$ as well.
Thus $\s \in \Int ^\CT _L (X)$ and (\ref{EQA020}) is true.

(b) First, using induction we prove the following auxiliary statement.
\begin{cla}\label{TA026}
For each  $L$-structure $\la X,\r \ra$ and each formula $\f \in \Form _{L_{\infty  \o }}$
we have
\begin{equation}\label{EQA026}
\forall \vec x\in {}^{\k }X \;\Big(    \la X,\r ^c \ra \models  \f  ^c [\vec x] \Leftrightarrow \la X,\r \ra \models  \f [\vec x]\Big).
\end{equation}
\end{cla}
\dok
Let $\vec x\in {}^{\k }X$. Then
$\la X,\r ^c \ra \models(v_\a =v_\b ) ^c [\vec x] $ iff $ \la X,\r ^c \ra \models ( v_\a =v_\b )[\vec x]$ iff $x_\a =x_\b$
iff $ \la X,\r  \ra \models ( v_\a =v_\b )[\vec x]$.
Also,
$\la X,\r ^c \ra \models(R_i (v_{\a _1}, \dots , v_{\a _{n_i}})) ^c [\vec x] $ iff $ \la X,\r ^c \ra \models  \neg  R_i (v_{\a _1}, \dots , v_{\a _{n_i}}) [\vec x]$
iff $\la x_{\a _1}, \dots , x_{\a _{n_i}}\ra \not\in \r ^c$ iff $\la x_{\a _1}, \dots , x_{\a _{n_i}}\ra \in \r $ iff
$\la X,\r  \ra \models R_i (v_{\a _1}, \dots , v_{\a _{n_i}}) [\vec x] $.

Suppose that (\ref{EQA026}) holds for a formula $\f$ and let $\vec x\in {}^{\k }X$. Then
$\la X,\r ^c \ra \models(\neg \f  ) ^c [\vec x] $
iff  not $ \la X,\r ^c \ra \models   \f  ^c   [\vec x]$, iff  not $ \la X,\r  \ra \models   \f    [\vec x]$, iff  $ \la X,\r  \ra \models (\neg  \f )   [\vec x]$.
$\la X,\r ^c \ra \models(\forall v_\a \; \f ) ^c [\vec x] $ iff $\la X,\r ^c \ra \models \forall v_\a \; \f  ^c [\vec x] $
iff for each $x\in X$ we have $\la X,\r ^c \ra \models \f  ^c [\vec x _{\la \a , x\ra}] $, that is, by (\ref{EQA026}), $\la X,\r  \ra \models \f  [\vec x _{\la \a , x\ra}] $,
iff $ \la X,\r \ra \models  (\forall v_\a \; \f )  [\vec x]$.
$\la X,\r ^c \ra \models(\exists v_\a \; \f ) ^c [\vec x] $ iff $\la X,\r ^c \ra \models \exists v_\a \; \f  ^c [\vec x] $
iff there is $x\in X$ such that $\la X,\r ^c \ra \models \f  ^c [\vec x _{\la \a , x\ra}] $, that is, by (\ref{EQA026}),
 $\la X,\r  \ra \models \f  [\vec x _{\la \a , x\ra}] $,
iff $ \la X,\r \ra \models  (\exists v_\a \; \f )  [\vec x]$.

Let $\Phi \subset \Form _{L_{\infty  \o }}$ and suppose that (\ref{EQA026}) holds for each formula $\f \in \Phi$.
$\la X,\r ^c \ra \models(\bigwedge \Phi  ) ^c [\vec x] $ iff $ \la X,\r ^c \ra \models \bigwedge \{\f  ^c :\f \in  \Phi \} [\vec x] $,
iff for each $\f \in  \Phi$ we have $ \la X,\r ^c \ra \models \f  ^c [\vec x] $, that is, by (\ref{EQA026}),  $ \la X,\r  \ra \models \f  [\vec x] $,
iff $\la X,\r \ra \models(\bigwedge \Phi  ) [\vec x] $.
$\la X,\r ^c \ra \models(\bigvee \Phi  ) ^c [\vec x] $ iff $ \la X,\r \ra \models  \bigvee \{\f  ^c :\f \in  \Phi \} [\vec x]$
iff for some $\f \in  \Phi$ we have $ \la X,\r ^c \ra \models \f  ^c [\vec x] $, that is, by (\ref{EQA026}),  $ \la X,\r  \ra \models \f  [\vec x] $,
iff $\la X,\r \ra \models(\bigvee \Phi  ) [\vec x] $.
\kdok
By Claim \ref{TA026}, the sets $(\Int ^\CT _L (X)) ^c = \{ \r ^c : \forall \f \in \CT  \; \la X,\r \ra\models \f\}$
and $\Int ^{\CT  ^c}_L (X) = \{ \r ^c : \forall \f \in \CT  \; \la X,\r ^c \ra\models \f  ^c \}$ are equal.

Statements (c) and (d) follow from (a), (b) and Theorem \ref{TA018}(a) and (c).
\hfill $\Box$
\begin{ex}\rm \label{EXA015}
Reversibility, complete theories and elementary equivalence.

If $\r \in \Int _L (X)$ and $\Th (\la X, \r \ra)$  is the corresponding first-order theory, then $[\r ]_\equiv :=\Int ^{\Th (\la X, \r \ra)}_L (X) $
is the set of interpretations $\s\in \Int _L (X)$ such that the structures $\la X, \r \ra$ and $\la X, \s \ra $ are elementarily equivalent.
We show that, regarding the relationship between the sets $[\r ]_\equiv $ and $\Rev _L (X)$, everything is possible.

1. If $\Q=\la Q,\r\ra$ is the rational line, then $\Th (\Q)$ is the theory of dense linear orders without end points
which is $\o$-categorical and, hence, $[\r ]_\equiv =[\r ]_{\cong}\subset \Rev _L (Q)$. By Fact \ref{T4312} $[\r ]_{\cong}$ is an antichain so
each element of the set $\Int ^{\Th (\Q )}_L (Q)$ is both a maximal and a minimal element of the set
$\Int ^{\Th (\Q )}_L (Q)$.

2. If $\BG =\la G,\r\ra$ is the countable universal homogeneous graph (also called the Rado graph, the Erd\H{o}s-R\'enyi graph \cite{Erdos2}), then the theory $\Th (\BG )$ is $\o$-categorical
and the structure $\BG $ is not reversible (since deleting of one of its edges produces an isomorphic copy of $\BG$, see \cite{Camer}).
Thus $\Int ^{\Th (\BG )}_L (G )\cap \Rev _L(G)=\emptyset$ and the set $\Int ^{\Th (\BG )}_L (G )$ has neither minimal
nor maximal elements.

3. It is well known that the theory $\CT$ of one equivalence relation having exactly one equivalence class of size $n$, for each $n\in \N$, is complete.
For a cardinal $\k \leq \o$ let $\E _\k=\la \o ,\r _\k\ra$ be a countable model of $\CT$ having exactly $\k$-many  infinite equivalence classes.
It is known that $\Int ^{\CT }_L (\o )=\bigcup _{\k \leq \o}[\r _\k]_{\cong }(=[\r _0]_\equiv )$ and, hence, $\CT$ is not an $\o$-categorical theory.
By \cite{KM3}, an equivalence relation  is reversible iff the number of equivalence classes of the same size
is finite or all equivalence classes are finite and their sizes form a reversible sequence. Thus the structures $\E _n$, $n<\o$, are reversible,
while $\E _\o$ is not (even weakly) reversible.
So we have $\Int ^{\CT }_L (\o ) \cap \Rev _L (\o )=\bigcup _{n\in \o}[\r _n ]_{\cong }$
and $\Int ^{\CT }_L (\o ) \setminus \Rev _L (\o )=[\r _\o ]_{\cong }$.

We show that $\Max ( \Int ^{\CT }_L (\o ))=[\r _0 ]_{\cong }\cup [\r _1 ]_{\cong }$.
Suppose that $n\in \{ 0,1 \}$ and $\r _n \varsubsetneq \s \in  \Int ^{\CT }_L (\o )  $. For $k\in \N$, let
$C_k$ be the unique equivalence class of size $k$ determined by $\r _n$. Since $\r _n \varsubsetneq \s$, the equivalence classes corresponding to $\s$
are unions of those corresponding to $\r$ and, in addition, there is the minimal
$k_0$ such that, in $\la \o ,\s \ra$, the class $C_{k_0}$ is joined with some another $\r$-class.
But then there is no $\s$-class of size $k_0$, which contradicts our assumption that  $\s \in  \Int ^{\CT }_L (\o )  $.
If $\k \geq 2$, then $\r _\k \subset \s$, for some $\s \cong \r _1 $ (we join $\k$ infinite classes into one).
Similarly we show that $\Min ( \Int ^{\CT }_L (\o ))=[\r _0 ]_{\cong }$ and that for $\k \geq 1$ there are no minimal elements of $\Int ^{\CT }_L (\o )$
below $\r _\k$ (split infinite $\r _\k$-classes into infinite parts).

Since $\r _m \preccurlyeq _c \r _n $, for $1\leq n\leq m \leq \o$,
the suborder $\{ [\r ]_{\sim _c}: \r \in \Int ^{\CT }_L (\o )\}=\{ [\r _n]_{\cong }:n\in \o \} \cup \{[\r _\o]_{\sim _c}\}$
of the condensation order $\la \Int _L (\o ) /\!\!\sim _c, \leq _c\ra$ is isomorphic to the disjoint union of the one element poset
(corresponding to $[\r _0]_{\cong }$) and the chain of the type $1+ \o ^*$, with the maximum $[\r _1]_{\cong }$ and minimum $[\r _\o]_{\sim _c}$.
\end{ex}
\section{Theories having extreme interpretations}\label{S4}
Example \ref{EXA015} shows that some sets of the form  $\Int ^{\CT }_L (X )$ have neither minimal nor maximal elements.
In this section we give some syntactical conditions providing extreme interpretations in that sense.
First, in order to provide maximal interpretations we define the class of {\it $R$-positive $L_{\infty  \o }$-formulas} by
$ \CP  :=\bigcup _{\xi \in \Ord} \CP  _\xi$,
where
\begin{eqnarray*}
 \CP  _0        &\! =  \!  & \At _{L } \cup \{ \neg v_\a =v_\b : \a ,\b \in \k  \},\\
 \CP  _{\xi +1} &\! =  \!  &  \CP  _\xi  \; \cup \; \{ \forall v_\a \; \f : \a\in \k \land  \f\in  \CP  _\xi  \}
                                        \; \cup \; \{ \exists v_\a \; \f : \a\in \k \land \f\in  \CP  _\xi  \}\\
                 &\! \cup \!& \textstyle \{ \bigwedge \Phi :  \Phi \subset  \CP  _\xi  \} \;  \cup  \; \{ \bigvee \Phi :  \Phi \subset  \CP  _\xi  \} ,\\
 \CP  _{\gamma} &\! =   \! & \textstyle  \bigcup _{\xi <\gamma } \CP _\xi , \mbox{ for a limit ordinal } \gamma ,
\end{eqnarray*}
\noindent
and  the class of $L_{\infty  \o }$-formulas $\CF :=\bigcup _{\xi \in \Ord}\CF _\xi$, where
\begin{eqnarray*}
\CF _0        &\! =  \!  &  \CP  \cup \{ \neg R_i (v_{\a _1}, \dots , v_{\a _{n_i}}) : i\in I \land  \la \a _1, \dots , \a _{n_i} \ra \in \k ^{n_i} \},\\
\CF _{\xi +1} &\! =  \!  & \CF _\xi  \; \cup \; \{ \forall v_\a \; \f : \a\in \k \land  \f\in \CF _\xi  \}   \\
              &\! \cup \!& \textstyle \{ \bigwedge \Phi :  \Phi \subset \CF _\xi  \} \;  \cup  \; \{ \bigvee \Phi :  \Phi \subset \CF _\xi \land |\Phi |<\o \} , \\
\CF _{\gamma} &\! =   \! & \textstyle  \bigcup _{\xi <\gamma }\CF _\xi , \mbox{ for a limit ordinal } \gamma .
\end{eqnarray*}
Concerning minimal interpretations, let us define
the class of {\it $R$-negative $L_{\infty  \o }$-formulas} by
$\CN =\bigcup _{\xi \in \Ord} \CN   _\xi$,
where
\begin{eqnarray*}
 \CN   _0        &\! =  \!  & \{ \neg R_i (v_{\a _1}, \dots , v_{\a _{n_i}}) : i\in I \land  \la \a _1, \dots , \a _{n_i} \ra \in \k ^{n_i} \}\\
                 &\! \cup \!& \{ v_\a =v_\b : \a ,\b \in \k  \} \cup  \{ \neg v_\a =v_\b : \a ,\b \in \k  \}       ,\\
 \CN   _{\xi +1} &\! =  \!  &  \CN   _\xi  \; \cup \; \{ \forall v_\a \; \f : \a\in \k \land  \f\in  \CN   _\xi  \}
                                        \; \cup \; \{ \exists v_\a \; \f : \a\in \k \land \f\in  \CN   _\xi  \}\\
                 &\! \cup \!& \textstyle \{ \bigwedge \Phi :  \Phi \subset  \CN   _\xi  \} \;  \cup  \; \{ \bigvee \Phi :  \Phi \subset  \CN   _\xi  \} ,\\
 \CN   _{\gamma} &\! =   \! & \textstyle  \bigcup _{\xi <\gamma } \CN  _\xi , \mbox{ for a limit ordinal } \gamma ,
\end{eqnarray*}
and let $\CG $ be the class of $L_{\infty  \o }$-formulas $\bigcup _{\xi \in \Ord}\CG _\xi$, where
\begin{eqnarray*}
\CG _0        &\! =  \!  &  \CN   \cup \{ R_i (v_{\a _1}, \dots , v_{\a _{n_i}}) : i\in I \land  \la \a _1, \dots , \a _{n_i} \ra \in \k ^{n_i} \},\\
\CG _{\xi +1} &\! =  \!  & \CG _\xi  \; \cup \; \{ \forall v_\a \; \f : \a\in \k \land  \f\in \CG _\xi  \}   \\
              &\! \cup \!& \textstyle \{ \bigwedge \Phi :  \Phi \subset \CG _\xi  \} \;  \cup  \; \{ \bigvee \Phi :  \Phi \subset \CG _\xi \land |\Phi |<\o \} , \\
\CG _{\gamma} &\! =   \! & \textstyle  \bigcup _{\xi <\gamma }\CG _\xi , \mbox{ for a limit ordinal } \gamma .
\end{eqnarray*}
\begin{te}\label{TA014}
Let $L$ be a relational language, $X$ a non-empty set and $\CT $ a set of $L_{\infty  \o}$-sentences
such that $\Int _L ^{\CT }(X)\neq\emptyset$. Then

(a) If $\CT \subset \CF $, then the set $\Int _L ^{\CT }(X)$ is union-complete and
$\Max ( \Int _L ^{\CT }(X))$ is its co-dense subset consisting of reversible interpretations;

(b) If $\CT \subset \CG $, then $\Int _L ^{\CT }(X)$ is intersection-complete and
$\Min (\Int _L ^{\CT }(X))$ is its dense subset consisting of reversible interpretations.
\end{te}
A proof is given in the sequel. First by induction we prove the following claim.
\begin{cla}\label{TA013}
(a) The formulas from the class $ \CP $ are preserved under condensations.

(b) For each formula $\f \in \CF$,  each chain $\L\subset \Int _L (X)$ and  valuation
$\vec x \in {}^{\k }X $ we have:
\begin{equation}\label{EQA014}\textstyle
\Big(\forall \r \in \L \;\;\la X,\r \ra\models \f [\vec x]\Big)
                               \Rightarrow \la X,\bigcup \L \ra\models \f [\vec x] .
\end{equation}
\end{cla}
\dok
(a) Let $\X$ and $\Y$ be $L$-structures and $f:X\rightarrow Y$ a condensation. By induction we show that (\ref{EQA021})
holds for each formula $\f\in  \CP $. First, clearly, homomorphisms preserve all atomic formulas.
If $\vec x\in {}^{\k}X$ and  $\X \models (\neg v_\a =v_\b)[\vec x ]$, that is $x _\a \neq x_\b$, then, since $f$ is an injection,
$f(x _\a) \neq f(x_\b )$, that is  $\Y \models (v_\a =v_\b)[f \vec x ]$.

Suppose that (\ref{EQA021}) holds for a formula $\f\in  \CP $; let $\vec x\in {}^{\k}X$.
If $\X \models (\forall v_\a \, \f )[\vec x ]$, then for each $x\in X$ we have $\X \models \f [\vec x _{\la \a , x \ra} ]$ and, by (\ref{EQA021}),
$\Y \models \f [(f \vec x) _{\la \a , f(x) \ra} ]$.  Since $f$ is a surjection, for each $y\in Y$ there is $x\in X$ such that $y=f(x)$ and, hence,
$\Y \models \f [(f \vec x) _{\la \a , y \ra} ]$. Thus $\Y \models (\forall v_\a \; \f )[f \vec x ]$.

If $\X \models (\exists v_\a \; \f )[\vec x ]$,
then for some $x\in X$ we have $\X \models \f [\vec x _{\la \a , x \ra} ]$
which  by (\ref{EQA021}) implies $\Y \models \f [(f \vec x) _{\la \a , f(x) \ra} ]$.
Thus, for $y=f(x)\in Y$ we have $\Y \models \f [(f \vec x) _{\la \a , y \ra} ]$ and, hence,  $\Y \models (\exists v_\a \; \f )[f \vec x ]$.

Let $\Phi \subset  \CP $, suppose that (\ref{EQA021}) holds for each formula $\f \in \Phi $ and let $\vec x\in {}^{\k}X$.
If $\X \models (\bigwedge \Phi )[\vec x ]$, then for each $\f \in \Phi $ we have
$\X \models \f [\vec x ]$ and, by (\ref{EQA021}),
$\Y \models \f [f \vec x ]$, which means that
$\Y \models (\bigwedge \Phi )[f \vec x ]$.
If $\X \models (\bigvee \Phi )[\vec x ]$, then for some $\f \in \Phi $ we have
$\X \models \f [\vec x ]$, that is
$\Y \models \f [f \vec x ]$, which implies that
$\Y \models (\bigvee \Phi )[f \vec x ]$.

(b) Let $\L$ be a non-empty chain in $\Int _L (X)$ and $\vec x \in {}^{\k }X $. (We note that then
$\bigcup \L =\la \bigcup _{\r \in \L}\r _i :i\in I\ra = \la \t _i :i\in I\ra =:\t \in \Int _L (X)$ and,
for an $i\in I$, the set $\L _i := \{ \r _i : \r \in \L \}$ is a chain in the algebra $\la \Int _{\la R _i\ra }(X),\subset\ra$.)

Let $\f \in \CP $ and suppose that $\la X,\r \ra\models \f [\vec x]$, for each $\r \in \L $.
If $\r \in \L$, then by Fact \ref{TA001}(a) the identity mapping
$\id _X :\la X, \r \ra\rightarrow \la X, \bigcup \L \ra$ is a condensation and, by (a),  preserves $\f$.
Thus, since $\la X,\r \ra\models \f [\vec x]$ we obtain $\la X,\bigcup \L \ra\models \f [\vec x]$.

If $\f := \neg R_i(v_{\a _1},\dots , v_{\a_{n_i}})$, then for any $\r \in \Int _L (X)$
we have  $\la X,\r \ra\models \f [\vec x]$ iff $\la x_{\a _1},\dots , x_{\a _{n_i}}\ra \not\in \r_i$.
Now, if $\la x_{\a _1},\dots , x_{\a _{n_i}}\ra \not\in \r_i$,
for each $\r \in \L $, then,
$\la x_{\a _1},\dots , x_{\a _{n_i}}\ra \not\in \t_i$, that is $\la X,\bigcup \L \ra\models \f [\vec x]$.

Suppose that the statement is true for a formula $\f \in \CF _\xi$.
Let $\L$ be a chain in $\Int _L (X)$ and $\vec x \in {}^{\k } X $.
If  for each $\r \in \L $ we have $\la X,\r \ra\models (\forall v_\a  \f) [\vec x]$,
that is  $\la X,\r \ra\models \f [\vec x _{\la \a , y\ra}]$, for all $y\in X$,
then for each $y\in X$ and each $\r \in \L $ we have $\la X,\r \ra\models \f [\vec x _{\la \a , y\ra}]$
and by the inductive hypothesis and (\ref{EQA014}) it follows that
$\la X,\bigcup \L \ra\models \f [\vec x _{\la \a , y\ra}]$.
This holds for all $y\in X$ so, $\la X,\bigcup \L \ra\models (\forall v_\a  \f) [\vec x ]$.

Let $\Phi \subset \CF _\xi$ and suppose that the statement is true for each formula $\f \in \Phi $.
Let $\L$ be a chain in $\Int _L (X)$ and $\vec x \in {}^{\k } X $.

If $\Phi =\{ \p _k : k\leq n \}$ and for each $\r \in \L $ we have $\la X,\r \ra\models (\bigvee _{k=1}^n \p _k )[\vec x]$,
then there are $k_0\leq n$ and a cofinal subset $\L _0$ of $\L$ such that $\la X,\r \ra\models \p _{k_0} [\vec x]$, for every $\r \in \L _0$
and, by the induction hypothesis, $\la X,\bigcup \L_0 \ra\models \p _{k_0} [\vec x]$.
By the cofinality of $\L _0$ we have $\bigcup \L_0=\bigcup \L$ and, hence, $\la X,\bigcup \L \ra\models (\bigvee _{k=1}^n \p _k ) [\vec x]$.

If for each $\r \in \L $ we have $\la X,\r \ra\models (\bigwedge \Phi ) [\vec x]$,
that is $\la X,\r \ra\models \f [\vec x]$, for all $\f\in \Phi$,
then, for each  $\f\in \Phi$ and $\r \in \L $ we have $\la X,\r \ra\models \f [\vec x]$,
so, by the induction hypothesis, $\la X,\bigcup \L \ra\models \f [\vec x]$.
Thus $\la X,\bigcup \L \ra\models (\bigwedge \Phi ) [\vec x]$.
\kdok

\noindent
{\bf Proof of Theorem \ref{TA014}(a)} Let $\L\subset \Int _L ^{\CT }(X)$ be a chain.
If $\f \in \CT$, then for each $\r \in \L$ we have $\r \in \Int _L ^{\CT }(X)$ and, hence, $\la X,\r \ra\models \f$, which, by (\ref{EQA014}),
 implies that $\la X,\bigcup \L \ra\models \f$.
So $\bigcup \L \in \Int _L ^{\CT }(X)$ and, thus, the set $\Int _L ^{\CT }(X)$ is union-complete.
The second statement follows from  Fact \ref{TA028} and Theorem \ref{TA027}(c).
\hfill $\Box$
\begin{cla}\label{TA036}
(a)  $\CN   = \{ \f  ^c : \f \in  \CP \}$;

(b) $\CG = \{ \f  ^c : \f \in \CF \}, \mbox{ up to logical equivalence}$.
\end{cla}
\dok
(a)
($\supset$) We show that for each $\xi \in \Ord$ and each $\f \in  \CP _\xi $ we have $\f  ^c \in  \CN   _\xi$. For $\xi =0$ we have:
$(v_\a =v_\b ) ^c := v_\a =v_\b \in   \CN  _0$, $(\neg v_\a =v_\b ) ^c := \neg (v_\a =v_\b) ^c := \neg v_\a =v_\b  \in   \CN  _0$,
and $(R_i (v_{\a _1}, \dots , v_{\a _{n_i}})) ^c := \neg  R_i (v_{\a _1}, \dots , v_{\a _{n_i}}) \in  \CN  _0$.

Suppose that the statement is true for all $\xi < \zeta$. If $\zeta $ is a limit ordinal, then, clearly, the statement is true for $\zeta$.
Let $\zeta =\xi +1$. If $\f \in  \CP _\xi$, then $\f  ^c \in  \CN  _\xi$ and, hence,
$(\forall v_\a \; \f ) ^c := \forall v_\a \; \f  ^c \in  \CN  _{\xi + 1} $ and
$(\exists v_\a \; \f ) ^c := \exists v_\a \; \f  ^c  \in  \CN  _{\xi + 1}$.

If $\Phi \subset  \CP _\xi$, then $\f  ^c \in  \CN  _\xi$, for all $\f \in \Phi$, and, hence, we have
$(\bigwedge \Phi  ) ^c :=\bigwedge \{\f  ^c :\f \in  \Phi \} \in  \CN  _{\xi + 1} $,
and
$(\bigvee \Phi  ) ^c := \bigvee \{\f  ^c :\f \in  \Phi \} \in  \CN  _{\xi + 1}$.

($\subset$) We show that for each $\xi \in \Ord$ and each $\p \in  \CN  _\xi $ there is $\f \in  \CP  _\xi$ such that $\p =\f  ^c$.
So $v_\a =v_\b$ is the formula $(v_\a =v_\b ) ^c $, $\neg v_\a =v_\b$ is the formula $( \neg v_\a =v_\b ) ^c $,
and $\neg  R_i (v_{\a _1}, \dots , v_{\a _{n_i}}) $ is the formula
$(R_i (v_{\a _1}, \dots , v_{\a _{n_i}})) ^c $.

Suppose that the statement is true for all $\xi < \zeta$. If $\zeta $ is a limit ordinal, then, clearly, the statement is true for $\zeta$.
Let $\zeta =\xi +1$. If $\p \in  \CN   _\xi$, then there is $\f \in  \CP  _\xi$ such that $\p =\f  ^c$.
Now
$\forall v_\a \; \p
= \forall v_\a \; \f  ^c
= (\forall v_\a \; \f ) ^c $, and $\forall v_\a \; \f \in  \CP  _{\xi +1}$.
Also
$\exists v_\a \; \p
= \exists v_\a \; \f  ^c
= (\exists v_\a \; \f ) ^c $,
and $\exists v_\a \; \f \in  \CP  _{\xi +1}$.

If $\Phi \subset  \CN   _\xi$, then for each $\p \in \Phi$ there is $\f_\p \in  \CP  _\xi$ such that $\p =\f _\p  ^c$.
So  $\bigwedge \{ \f_\p : \p \in \Phi \}\in  \CP  _{\xi +1}$ and
$(\bigwedge \{ \f_\p : \p \in \Phi \}) ^c = \bigwedge \{ \f_\p  ^c: \p \in \Phi \}= \bigwedge \Phi$. Also
$\bigvee \{ \f_\p : \p \in \Phi \}\in  \CP  _{\xi +1}$ and $(\bigvee \{ \f_\p : \p \in \Phi \}) ^c = \bigvee \{ \f_\p  ^c: \p \in \Phi \}= \bigvee \Phi$.

(b) ($\supset$) We show that for each $\xi \in \Ord$ and each $\f \in \CF _\xi $ we have $\f  ^c \in \CG _\xi$.
For $\xi =0$, if $\f \in  \CP $, then, by (a), $\f  ^c \in  \CN   \subset \CG _0$ and
$(\neg R_i (v_{\a _1}, \dots , v_{\a _{n_i}})) ^c $ is the formula $\neg \neg  R_i (v_{\a _1}, \dots , v_{\a _{n_i}}) $, which is equivalent to
$R_i (v_{\a _1}, \dots , v_{\a _{n_i}}) \in \CG _0$.

Suppose that the statement is true for all $\xi < \zeta$. If $\zeta $ is a limit ordinal, then, clearly, the statement is true for $\zeta$.
Let $\zeta =\xi +1$. If $\f \in \CF _\xi $, then $\f  ^c \in \CG _\xi$
and, hence, $(\forall v_\a \; \f ) ^c := \forall v_\a \; \f  ^c \in \CG _{\xi +1}$.

If $\Phi \subset \CF _\xi$, then $\f  ^c \in \CG _\xi$, for all $\f \in \Phi$,
and, hence, we have
$(\bigwedge \Phi  ) ^c :=\bigwedge \{\f  ^c :\f \in  \Phi \} \in \CG _{\xi + 1} $,
and
$(\bigvee \Phi  ) ^c := \bigvee \{\f  ^c :\f \in  \Phi \} \in \CG _{\xi + 1}$, if $|\Phi|<\o$.

($\subset$) We show that for each $\xi \in \Ord$ and each $\p \in \CG _\xi $ there is $\f \in \CF _\xi$ such that $\p =\f  ^c$.
For $\xi =0$, if $\p \in  \CN  $, then we apply (a).
Also, $ R_i (v_{\a _1}, \dots , v_{\a _{n_i}}) $ is equivalent to the formula
$(\neg R_i (v_{\a _1}, \dots , v_{\a _{n_i}})) ^c $.

Suppose that the statement is true for all $\xi < \zeta$. If $\zeta $ is a limit ordinal, then, clearly, the statement is true for $\zeta$.
Let $\zeta =\xi +1$.
If $\p \in \CG _\xi$, then there is $\f \in \CF _\xi$ such that $\p =\f  ^c$.
Now
$\forall v_\a \; \p
= \forall v_\a \; \f  ^c
= (\forall v_\a \; \f ) ^c $, and $\forall v_\a \; \f \in \CF  _{\xi +1}$.

If $\Phi \subset \CG _\xi$, then for each $\p \in \Phi$ there is $\f_\p \in \CF _\xi$ such that $\p =\f _\p  ^c$.
So  $\bigwedge \{ \f_\p : \p \in \Phi \}\in \CF  _{\xi +1}$ and
$(\bigwedge \{ \f_\p : \p \in \Phi \}) ^c = \bigwedge \{ \f_\p  ^c: \p \in \Phi \}= \bigwedge \Phi$.
If $|\Phi|<\o$, then $\bigvee \{ \f_\p : \p \in \Phi \}\in \CF  _{\xi +1}$ and
$(\bigvee \{ \f_\p : \p \in \Phi \}) ^c 
= \bigvee \Phi$.
\kdok

\noindent
{\bf Proof of Theorem \ref{TA014}(b)}
If $\L$ is a chain in the poset $\la \Int _L ^{\CT }(X), \subset \ra$, then, by Theorem \ref{TA027}(b) we have
$\L  ^c =\{ \r ^c : \r \in \L \} \subset \{ \r ^c : \r\in \Int _L ^{\CT }(X)\} =:(\Int _L ^{\CT }(X)) ^c =\Int _L ^{\CT  ^c}(X)$.
Since $\CT \subset \CG $, by Claim \ref{TA036}(b)
w.l.o.g.\ we assume that $\CT \subset \{ \f  ^c :\f \in \CF\}$ and, hence,  $\CT  ^c \subset \{ (\f ^c)^c :\f \in \CF\}$.
By Claim \ref{TA026}, for each  interpretation $\r \in \Int _L (X)$ and each $L_{\infty  \o }$-sentence $\f $
we have: $\la X,\r \ra \models  \f $  iff $\la X,\r  \ra \models  (\f ^c)^c $ so, w.l.o.g.\ again, we suppose that
$\CT  ^c \subset \CF$. Now, clearly, $\L  ^c$ is a chain in the poset $\la \Int _L ^{\CT  ^c}(X), \subset \ra$ and, by Theorem \ref{TA014}(a),
$\bigcup \L  ^c =\bigcup _{\r\in \L}\r ^c =(\bigcap _{\r\in \L}\r )^c \in \Int _L ^{\CT  ^c}(X)$ and, by Theorem \ref{TA027},
$\bigcap _{\r\in \L}\r =\bigcap \L\in \Int _L ^{\CT }(X)$.
The second statement follows from Fact \ref{TA028} and Theorem \ref{TA027}(c).
\hfill $\Box$
\begin{ex}\label{EXA016}\rm
Extreme partial orders.
Clearly, for the set of axioms of the theory of strict partial orders $\CT _{poset}=\{ \f _{irr}, \f _{tr}\}\subset \Sent _{L_b} $,
where $\f _{irr}:= \forall v_0 \, \neg R(v_0,v_0)$ and
$\f _{tr}:= \forall v_0 ,v_1, v_2 ( \neg R(v_0,v_1) \lor \neg R(v_1,v_2)\lor R(v_0,v_2))$ we have $\CT _{poset}\subset \CF \cap \CG$ and, hence,
the poset  $\P :=\la \Int _{L_b}^{\CT _{poset}}(X),\subset \ra$ of all strict partial orders on $X$ has all the properties from (a) and (b) of Theorem \ref{TA014}.
It is evident that $\Min \P =\{ \emptyset \}$ and, by Example \ref{EXA012}, this {\it antichain order} is the unique strongly reversible strict partial order on $X$.

The maximal elements of the poset $\P$ are exactly the {\it strict linear orders}. Namely, it is clear that
if $\la X,\r\ra$ is a strict linear order and $\r \varsubsetneq \r '$, then $\r '$ is not a strict partial order. On the other hand, by the Order
extension principle (i.e.\ the Szpilrajn extension theorem \cite{Szp}, following from Zorn's lemma), if $\r$ is a strict partial order on $X$, then there is a strict linear order $\r '$ on $X$ such that $\r' \supset \r $.

We note that, by a well-known theorem of Dushnik and Miller \cite{DM},
the poset $\P$ has the following property: each interpretation  $\r \in\Int _{L_b}^{\CT _{poset}}(X)$
is the intersection of a family of maximal elements of $\P$ and the minimal size of such a family is called the {\it Dushnik-Miller dimension} of the poset
$\la X,\r\ra$.
In \cite{DM} a poset is called reversible iff it is of dimension $\leq 2$, but it is easy to check that
the poset $\X =\la Z, < \ra$, where $Z$ is the set of integers and $<:=\{ \la 2n-1, 2n \ra :n\in \N\}$ is of dimension 2,
but not reversible in our sense. In \cite{Kuk} Kukiela has shown that Boolean lattices are reversible posets (in our sense), but, clearly,
lot of them have dimension $>2$.
\end{ex}
\begin{ex}\rm \label{EXA018}
The poset of interpretations of countable connected graphs
is union-complete but not intersection-complete, although the minimal elements are dense in it.
For the set of axioms of graph theory $\CT _{graph}=\{ \f _{irr}, \f _{sym}\}$, where $\f _{irr}:= \forall v_0 \, \neg R(v_0,v_0)$
and $\f _{sym}:= \forall v_0 ,v_1 (\neg R(v_0,v_1) \lor R(v_1,v_0))$ we have
$\CT _{graph}\subset\CF $ and
the $L_{\infty \o }$-sentence $\f_{conn}$ given by
$$\textstyle
\forall u,v \;\Big( u=v \lor \bigvee _{n\geq 2} \exists v_1, \dots ,v_n \;( u=v_1 \land v=v_n \land \bigwedge _{k=1}^{n-1}R( v_k,v_{k+1})\Big)
$$
and expressing that a graph is connected belongs to $ \CP $.
So $\CT _{graph}\cup \{ \f_{conn}\}\subset \CF $ and, by Theorem \ref{TA014}(a),
the poset  $\la \Int _{L_b}^{\CT _{graph}\cup \{ \f_{conn}\}}(\o), \subset \ra$
is union-complete. Since a graph is a tree iff it is a minimal connected graph, the minimal elements of our poset are exactly
the tree graph relations on $X$. Since every connected graph admits a spanning tree (it is an easy application of Zorn's lemma; see \cite{Souk}), our poset has dense set of minimal elements.
For $k\in \o$, let $\BG _k =\la \o \cup \{ \o\}, \r _k\ra$, where $\r _k= \{ \{ n,n+1\} : n\in \o\}\cup \{ \{ n,\o  \} :n\geq k  \}$.
It is evident that the graphs $\BG _k$ are connected and $\r _0 \varsupsetneq \r _1 \varsupsetneq \r _2 \varsupsetneq \dots $, but the graph
$\BG _\o =\la \o \cup \{ \o\}, \bigcap _{k\in \o }\r _k\ra$ is disconnected and, hence, the poset is not intersection-complete.
\end{ex}
\section{Omitting finite substructures}\label{S6}
A class $\K \subset \Mod _L$ is called a {\it universal class} iff it is axiomatizable by a finite set of universal ($\Pi ^0_1$) sentences
iff there exists a finite set of finite $L$-structures $\{\F _k : k\leq n\}\subset \Mod _L$ such that $\X \in \K$ iff $\F_k\not\hookrightarrow \X$, for all
$k\leq n$ (see \cite{Tar,Vau,Fra1,Fra}). Here, using that concept, we show that forbidding finite structures provides a large zoo of reversible structures.
\begin{fac}\label{TA015}
For each finite $L$-structure $\F $ there is an $L_{\infty  \o }$-sentence $\p _{\F \hookrightarrow}$
such that for each $L$-structure $\Y$  we have:
$\F \hookrightarrow \Y $ iff $\;\Y \models \p _{\F \hookrightarrow}$.
If, in addition, the language $L$ is finite, then the sentence $\neg \p _{\F \hookrightarrow}$ is logically equivalent to a $\Pi ^0_1$ sentence
$\eta _{\F \not\hookrightarrow}$.
\end{fac}
\dok
Let $L=\la R_i: i\in I\ra$, where $\ar (R_i )=n_i$, for $i\in I$, and w.l.o.g.\ suppose that $\F =\la m, \la R_i ^\F :i\in I\ra \ra \in \Mod _L$,
where $m=\{0,\dots ,m-1\}\in \N$. Let $\chi _{R_i ^\F }:m ^{n_i}  \rightarrow 2$, $i\in I$, be
the characteristic functions of the sets $R_i ^\F \subset m ^{n_i}$ and let
$\f _\F (v_0, \dots , v_{m-1})$ be the $L_{\infty \o }$-formula defined by
\begin{equation}\label{EQA017}\textstyle
\f _\F (\bar v ):=
\bigwedge _{0\leq j<k <m}v_j\neq v_k \;\;\land \;\;\bigwedge _{i\in I}
\bigwedge _{\bar x \in m^{n_i}}R_i(v_{x_0}, \dots , v_{x_{n_i -1}})^{\chi _{R_i ^\F }(\bar x)},
\end{equation}
where, by definition, $\eta ^1:=\eta$ and $\eta ^0:=\neg \eta$.
We show first that
\begin{equation}\label{EQA018}
\F \models \f _\F [0,1,\dots ,m-1].
\end{equation}
For $j<m$, under the valuation $\la 0,1,\dots ,m-1 \ra$ the variable $v_j$ obtains the value $j$ and, hence,
$\F \models (\bigwedge _{0\leq j<k <m}v_j\neq v_k )\; [0,1,\dots ,m-1]$ is true.
Let $i\in I$ and $\bar x =\la x_0 , \dots ,x_{n_i -1}\ra\in m^{n_i}$. Then
$\F \models R_i(v_{x_0}, \dots , v_{x_{n_i -1}})^{\chi _{R_i ^\F }(\bar x)}\; [0,1,\dots ,m-1]$ iff
$\F \models R_i [x_0, \dots , x_{n_i -1}]^{\chi _{R_i ^\F }(\bar x)}$ iff
$(\chi _{R_i ^\F }(\bar x)=1 \land \bar x \in R_i^{\F }) \lor (\chi _{R_i ^\F }(\bar x)=0 \land \bar x \not\in R_i^{\F })$ which is true.
Thus (\ref{EQA018}) is proved.

Let $\p _{\F \hookrightarrow}:=\exists \bar v \; \f _\F (\bar v)$.
If $\Y\in \Mod _L$ and $f:\F \hookrightarrow \Y$, then by (\ref{EQA018})
we have $\Y \models \f _\F [f(0),\dots ,f(m-1)]$ and, since
$\bar y:= \la f(0),\dots ,f(m-1)\ra \in Y^m$ we have $\Y \models \exists \bar v \; \f _\F (\bar v)$, that is, $\Y \models \p _{\F \hookrightarrow}$.
Conversely, let $\bar y =\la y_0,\dots y_{m-1}\ra\in Y^m$ and $\Y \models \f _\F [\bar y]$. Since under the valuation $\bar y$ the variable $v_j$ obtains the value $y_j$,
by (\ref{EQA017}) $y_0,\dots y_{m-1}$ are different elements of $Y$ and, hence, the mapping $f:m \rightarrow Y$ defined by $f(j)=y_j$, for $j<m$,
is an injection. For a proof that $f:\F \rightarrow \Y$ is a strong homomorphism we take $i\in I$ and $\bar x :=\la j_0,\dots ,j_{n_i -1}\ra\in m ^{n_i}$ and show that
$$
\la j_0,\dots ,j_{n_i -1}\ra\in R^\F _i \Leftrightarrow \la y_{j_0},\dots ,y_{j_{n_i -1}}\ra\in R^\Y _i .
$$
Since $\Y \models \f _\F [\bar y]$, by (\ref{EQA017}) for $\bar x$ we have
$\Y \models R_i(v_{j_0}, \dots , v_{j_{n_i -1}})^{\chi _{R_i ^\F }(\bar x)} [\bar y]$, that is
$\Y \models R_i[y_{j_0}, \dots , y_{j_{n_i -1}}]^{\chi _{R_i ^\F }(\la j_0,\dots ,j_{n_i -1}\ra)} $, thus
$\la y_{j_0},\dots ,y_{j_{n_i -1}}\ra\in R^\Y _i $ if and only if $\chi _{R_i ^\F }(\la j_0,\dots ,j_{n_i -1}\ra)=1$
iff $\la j_0,\dots ,j_{n_i -1}\ra\in R^\F _i$ and that's it.

If $|L|<\o$, then the sentence $\neg \p _{\F \hookrightarrow}$ is equivalent to the $\Pi ^0_1$ sentence
$\textstyle
\eta _{\F \not\hookrightarrow }:=
\forall \bar v\;
(\bigvee _{0\leq j<k <m}v_j = v_k \;
\lor \;\bigvee _{i\in I} \bigvee _{\bar x \in m^{n_i}} R_i(v_{x_0}, \dots , v_{x_{n_i -1}})^{1-\chi _{R_i ^\F }(\bar x)})
$.
\hfill $\Box$
\begin{te}\label{TA029}
Let $L$ be a finite language, $\CT$ an $L_{\infty  \o }$-theory and $\F _j$, $j\in\! J $, finite $L$-structures such that the poset
$\P := \la \Int _L^{\CT \cup \{ \eta _{\F _j \not\hookrightarrow}: j\in J\}}(X),\subset \ra$ is non-empty.

(a) If $\CT\subset \CF$, then the poset $\P$ is union-complete and $\Max \P $ is a co-dense set in $\P$ consisting of reversible interpretations;

(b) If $\CT\subset \CG$, then the poset $\P$ is intersection-complete and $\Min \P $ is a dense set in $\P$ consisting of reversible interpretations;

(c) $\t \in \Max  (\Int _L^{\CT \cup \{ \eta _{\F _j \not\hookrightarrow}: j\in J\}}(X))\;$  iff 
    $\;\t ^c \in \Min ( \Int _L^{\CT ^c \cup \{ \eta _{\F^c _j \not\hookrightarrow}: j\in J\}}(X))$.
\end{te}
\dok
Since $ \Pi ^0_1\subset \CF \cap \CG$, (a) and (b) follow from Fact \ref{TA015} and Theorem \ref{TA014}.

(c) It is easy to check that for each $p\in \{ 0,1\}$ we have $(R_i (\bar v)^{1-p})^c \leftrightarrow R_i (\bar v)^{p}$ and, also,
that $\chi _{R_i ^\F }(\bar x) =1-\chi _{R_i ^{\F ^c} }(\bar x)$, which implies that
$(\eta _{\F \not\hookrightarrow })^c\leftrightarrow \eta _{\F ^c \not\hookrightarrow }$. So the statement
follows from Theorem \ref{TA027}(d).
\hfill $\Box$
\paragraph{Maximal $\BK _n$-free graphs}
In the sequel,  for convenience, for a graph $\X =\la X , \rho \ra$ the relation $\rho $ will be identified with the corresponding
set of two-element subsets of $X$, $\{\{ x,y \}\in [X]^2 : \la x,y \ra \in \rho \}$ and $\X^{gc}$ will denote the graph-complement, $\la X , [X]^2\setminus \rho \ra$,
of the graph $\X $. For a set $Y\subset X$, the subgraph $\la Y,\rho \upharpoonright Y\ra$ of $\X$ will be sometimes denoted by $Y$.
For a cardinal $\nu$,  $\BK _\nu$ will denote the {\it complete graph} of size $\nu$, and $\E_\nu$  the graph with $\nu$ vertices and no edges. Clearly, $\E_\nu=\BK_\nu^{gc}$.

If $\F$ is a finite graph which is not complete, then, trivially, $X^2 \setminus \Delta _X$
is the unique maximal element of the poset $\Int _{L_b}^{\CT _{graph}\cup \{ \eta _{\F \not\hookrightarrow}\}}(X)$ and here we consider what forbidding $\BK _n$'s produce.   By  Theorem \ref{TA029},
the poset $\Int _{L_b}^{\CT _{graph}\cup \{ \eta _{\BK _n \not\hookrightarrow}\}}(X)$ has maximal elements, they are reversible and, clearly,  different from
$X^2 \setminus \Delta _X$.
We recall
that a graph is called {\it $\BK _n$-free} iff it has no subgraphs isomorphic to $\BK _n$; trivially, the graphs $\BK _m$, $m<n$, are maximal $\BK _n$-free graphs.
\begin{cla}\label{TA040}
Let $n\geq 3$ and let $\X =\la X , \rho \ra$ be a $\BK _n$-free graph. Then

(a) $\X $ is a maximal $\BK _n$-free graph iff
\begin{equation}\label{EQA090}
\forall \{ x,y \}\in [X]^2 \setminus \rho \;\; \exists K\in [X]^n \;\; [K]^2 \setminus \rho =\{\{ x,y \}\}.
\end{equation}

(b) If $\X $ is a maximal $\BK _n$-free graph and $|X|\geq n-1$, then
\begin{equation}\label{EQA091}
\forall x\in X \;\; \exists K\in [X\setminus \{ x \}]^{n-2} \;\;\; \{x\}\cup K \cong \BK _{n-1}.
\end{equation}

(c) If $\X $ is a maximal $\BK _n$-free graph, $|X|\geq n-1$, $\{ Y_x :x\in X\}$ is a family of non-empty sets,
$Y:=\bigcup _{x\in X}\{ x \}\times Y_x$ and
\begin{equation}\label{EQA092}
\sigma =\Big\{ \{ \la x, y \ra , \la x', y' \ra\}\in [Y]^2 : \{ x,x' \}\in \rho \Big\},
\end{equation}
then $\Y =\la Y , \sigma \ra$ is a maximal $\BK _n$-free graph.

(d) $\X $ is a maximal $\BK _n$-free graph iff
 $\X ^c$ is a minimal $\la n,\Delta _n\ra$-free reflexive graph iff
$\X ^{gc}$ is a minimal $\E _n$-free graph.
\end{cla}

\dok
(a)  If $|X|<n$, then (\ref{EQA090}) holds iff $\rho =[X]^2$ iff $\la X, \rho \ra\cong \BK _{|X|}$. Let $|X|\geq n$.

If $\X$ is maximal and $\{x,y\}\in[X]^2\setminus\rho $, then the graph
$\la X,\rho\cup\{\{x,y\}\}\ra$ is not $\BK_n$-free, which means that there is a set $K\in[X]^n$
such that $x,y\in K$ and
$\la K,(\rho \cup\{\{x,y\}\})\upharpoonright K\ra\cong\BK_n$, which implies that $[K]^2\setminus\rho=\{\{x,y\}\}$.

Conversely, if (\ref{EQA090}) holds, then for any $\{x,y\}\in[X]^2\setminus\rho $ there is $K\in[X]^n$ such  that
$\la K,(\rho \cup\{\{x,y\}\})\upharpoonright K\ra\cong\BK_n$, thus $\X$ is a maximal $\BK_n$-free graph.

(b) If $|X|=n-1$, then $\X\cong \BK _{n-1}$ and (\ref{EQA091}) is evident.
Let $|X|\geq n$ and $x\in X$. If $\{ x,y \}\not\in \rho $, for some $y \in X\setminus \{ x \}$, then by (\ref{EQA090})
there is a set $K'=\{ x,y,x_1,\dots ,x_{n-2} \}\in [X]^n $ such that $[K']^2 \setminus \rho =\{\{ x,y \}\}$ and for
$K:=\{ x_1,\dots ,x_{n-2} \}\in [X\setminus \{ x \}]^{n-2} $ we have $\{x\}\cup K \cong \BK _{n-1}$.

If $\{ x,y \}\in \rho $, for all $y \in X\setminus \{ x \}$, then, since $|X|\geq n$, there is a pair $\{ u,v \}\in [X\setminus \{ x \}]^2 \setminus \rho $
and, by (\ref{EQA090}), there is a set  $K=\{ u,v, x_1,\dots ,x_{n-2}\}\in [X]^n $ such that $[K]^2 \setminus \rho =\{\{ u,v \}\}$.
Now, if $x\not\in \{ x_1,\dots ,x_{n-2}\}$, then $\{ x\} \cup \{ x_1,\dots ,x_{n-2}\}\cong \BK _{n-1}$ and (\ref{EQA091}) is true.
If $x=x_j$, for some $j\leq n-2$, then $\{ x\} \cup \{u, x_1,\dots ,x_{j-1}, x_{j+1}, \dots ,x_{n-2}\}\cong \BK _{n-1}$ and (\ref{EQA091}) is true again.

(c) Suppose that $\{ \la x_i ,y_i \ra : 1\leq i \leq n\}$ is a copy of $\BK _n$ in $\Y$; then, by (\ref{EQA092}),
$\{  x_i : 1\leq i \leq n\}$ would be a copy of $\BK _n$ in $\X$, which contradicts our assumption. Thus $\Y $ is a $\BK _n$-free graph.

Suppose that $\la Y ,\tau \ra$ is a $\BK _n$-free graph, where $\sigma \varsubsetneq \tau$.
Let $\{\la x, y \ra , \la x', y' \ra\}\in \tau\setminus \sigma$ . If $x=x'$, then by (b) there is a set
$K=\{ x_1,\dots ,x_{n-2}\}\in [X\setminus \{ x \}]^{n-2}$ such that
$\{x\}\cup K \cong \BK _{n-1}$. For $j\leq n-2$ we choose $y_i\in Y_{x_i}$ and,  by (\ref{EQA092}),
$\{\la x, y \ra , \la x, y' \ra\} \cup \{ \la x_j ,y_j  \ra :j\leq n-2\}$ is a copy of $\BK _n$ in $\la Y ,\tau \ra$, contrary to our assumption.

If $x\neq x'$, then $\{ x,x' \}\in [X]^2 \setminus \rho $ and by (a) there is $K=\{ x,x', x_1,\dots ,x_{n-2} \}$ $\in [X]^n$ such that $[K]^2 \setminus \rho =\{\{ x,x' \}\}$. Again, for $j\leq n-2$ we choose $y_i\in Y_{x_i}$ and,  by (\ref{EQA092}),
$\{\la x, y \ra , \la x', y' \ra\} \cup \{ \la x_j ,y_j  \ra :j\leq n-2\}$ is a copy of $\BK _n$ in $\la Y ,\tau \ra$, contrary to our assumption.
Thus $\Y $ is a maximal $\BK _n$-free graph.

(d) Clearly, up to logical equivalence we have $\CT ^c_{graph}=\{ \varphi _{refl},\varphi _{sym}\}$ and $\BK _n ^c \cong\la n, \Delta _n\ra$.
Now the first claim follows from Theorem \ref{TA029}(c) and the second claim follows from the first one.
\kdok
\begin{ex}\label{EXA030}\rm
Claim \ref{TA040} provides a large jungle of extreme and, hence, reversible structures. So if $n\geq 3$, $\X\cong \BK _{n-1}$ and
$\{ Y_x :x\in X\}$ is family of non-empty sets, then the graph $\Y$ defined in Claim \ref{TA040}(c)  is a maximal $\BK _n$-free graph.
The reader will notice that $\Y$ is in fact the {\it complete $(n-1)$-partite graph} and that $\Y ^{gc}$ is a disjoint union of $n-1$ complete graphs, which is a minimal $\E _n$-free graph. If $|Y _x|=\omega $, for all $x\in X$, then $\Y ^{gc}$ is a reversible countable ultrahomogeneous graph from the list of Lachlan and Woodrow
(see Remark \ref{RA010}).

For $n=3$, the complete bi-partite graphs $\BK_{\nu,\omega }$, $\nu\leq\omega $, are maximal countable triangle-free graphs. In particular, the {\it star graph} $\S _\omega :=\BK_{1,\omega }$, is a maximal triangle-free graph. Furthermore, some maximal triangle-free graphs are not bi-partite,
for example the cycle graph $\C _5$. Also, by taking $\X\cong \C _5$ in Claim \ref{TA040}(c)
we obtain infinite maximal $\BK _3$-free graphs which are not bi-partite.

Of course, there are reversible $\BK _3$-free graphs which are not maximal $\BK _3$-free. For example, the {\it linear graph} $\BG _\omega =\la \omega ,\tau \ra$, where $\tau =\{ \{ n, n+1\} : n\in \omega \}$, is reversible since deleting an edge produces a disconnected graph.
\end{ex}
\paragraph{Maximal $\BK _n$-free graphs with all vertices of infinite degree}
In the context of graph theory, the sentence
$$\textstyle
\varphi _\infty := \forall v \;\; \bigwedge _{n\in \N} \exists v_1 , \dots , v_n \;
\Big(\bigwedge _{1\leq i < j \leq n}v_i \neq v_j \land \bigwedge _{1\leq i \leq n} R(v, v_i)\Big)
$$
says that each vertex of a graph has infinitely many neighbors. Since $\varphi _\infty \in \CP $,
by  Theorem  \ref{TA029} the poset $\Int _{L_b}^{\CT _{graph}\cup \{  \varphi _\infty, \eta _{\BK _n \not\hookrightarrow} \}}(X)$
has a co-dense set of maximal elements and they are reversible. Some such interpretations are given in Example \ref{EXA030}.
\begin{ex}\rm \label{EXA019}
The Henson graph $\BH _n$ is a maximal $\BK _n$-free graph with all vertices of infinite degree.
For $n\geq 3$, $\BH _n$ denotes the unique countable homogeneous universal $\BK _n$-free graph (the {\it Henson graph}, see \cite{henson}).
In order to recall a convenient characterization of $\BH _n$
we introduce the following notation:
if $\BG=\la G,\rho \ra$ is a graph and $n\geq 3$ let
$C_n(\BG ):=\{\la H,K \ra : K\subset H\in [G]^{<\omega } \land K \mbox{ is } \BK _{n-1}\mbox{-free}\}$ and
for $\la H,K \ra \in C_n(\BG )$ let
$$G^H_K :=\{v\in G\setminus H: \forall k\in K\; \{ v ,k \} \in \rho \land \forall h\in H\setminus K \; \{v, h \}  \notin \rho \}.$$
Now, by \cite{henson} we have: a countable graph $\BG=\la G,\rho \ra$ is isomorphic to $\BH _n$ iff
$\BG$ is $\BK_n$-free and $G^H_K\neq \emptyset$, for each $\la H,K \ra\in C_n(\BG )$.

We show that the Henson graph $\BH _n =\la G,\rho \ra$ is a maximal $\BK _n$-free graph.
Suppose that $\la G ,\rho '\ra $ is a $\BK_n$-free graph, where $\rho \varsubsetneq \rho '$ and $\{ a_1,a_2 \}\in \rho '\setminus \rho $.
By recursion we construct different elements $a_3, \dots, a_n \in G\setminus \{ a_1,a_2 \}$ such that
\begin{equation}\label{EQA029}
\forall k\in \{ 3,\dots ,n\} \;\;\forall i\in \{1,2,\dots ,k-1\} \;\;\{ a_i,a_k \}\in \rho .
\end{equation}
Let $k\in \{ 3,\dots ,n\}$ and suppose that the sequence $ a_1,a_2,\dots , a_{k-1} $ satisfies (\ref{EQA029}).
Then, since $\{ a_1,a_2 \}\not\in \rho $, for $H=K:=\{ a_1,a_2,\dots , a_{k-1} \}$ we have
$\BK _{n-1} \not\hookrightarrow \la K, \rho \upharpoonright K\ra$
and, hence, $\la H,K\ra\in C_n (\BH _n)$
so, by the characterization mentioned above,
there is $a_k\in G\setminus \{ a_1,a_2,\dots , a_{k-1} \}$ such that $\{ a_i,a_k \}\in \rho $, for all $i<k$.
So, the sequence $ a_1,a_2,\dots , a_k$ satisfies (\ref{EQA029}) and the recursion works.
But, since $\{ a_1,a_2 \}\in \rho '$ the vertices $a_1, \dots, a_n$ determine a subgraph of the graph $\la G ,\rho '\ra $ isomorphic to $\BK_n$,
which contradicts our assumption.

Since the star graph $\S _\omega $ (see Example \ref{EXA030}) is $\BK _n$-free, by the universality of $\BH _n$ there is a copy of $\S _\omega $ in $\BH _n$ and, hence,
$\BH _n$ contains a vertex of infinite degree. By the homogeneity of $\BH _n$ all vertices of $\BH _n$ are of infinite degree.
\end{ex}
\begin{rem}\label{RA010}\rm
By a well-known characterization of Lachlan and Woodrow \cite{Lach},
each countable ultrahomogeneous graph is isomorphic to one of the following:

- $\BG _{\mu \nu }$, the union of $\mu$ disjoint copies of $\BK _\nu$, where $\mu \nu=\omega $, - reversible iff $\mu <\omega $ or $\nu <\omega $ \cite{KM3};

- $\BG_{\Rado}$, the Rado graph - non-reversible (see Example \ref{EXA015});

- $\BH _n$, the Henson graph, for $n\geq 3$ - reversible (see Example \ref{EXA019});

- the graph-complements of these graphs - a graph is reversible iff its graph-complement is reversible (it is an easy consequence of Fact \ref{T4312}).
\end{rem}

\paragraph{Omitting extreme finite structures}
Clearly, the minimal elements of the set $\Int _L^{\{ \eta _{\F \not\hookrightarrow}\}}(X)$, in the sequel denoted shortly by
$\Int _L^{ \eta _{\F \not\hookrightarrow} }(X)$, will be different from the trivial one, $\la \emptyset :i\in I\ra $,
iff the forbidden structure $\F$ is minimal, that is isomorphic to $\la m, \la \emptyset:i\in I\ra\ra$, for some $m\in \N$.
Dually, $\Max ( \Int _L^{ \eta _{\F \not\hookrightarrow}}(X))\neq \{  \la X, \la X^{n_i}:i\in I\ra\ra\}$ iff
$\F \cong\la m, \la m^{n_i}:i\in I\ra\ra$. We give some examples of such restrictions.
\begin{cla}\label{TA041}
Let $m,n\in \N$, $L_n=\la R\ra$, where $\ar (R)=n$. Then

(a) If $\rho \in \Int _L^{ \eta _{\la m, \emptyset\ra \not\hookrightarrow}}(X) $,
then $\rho \in \Min (\Int _{L_n}^{\eta _{\la m, \emptyset\ra \not\hookrightarrow}}(X) )$ iff
\begin{equation}\label{EQA093}
\forall \bar x \in \rho \;\;\exists K \in [X]^m \;\;\rho \cap K^n =\{ \bar x\};
\end{equation}

(b) If $\rho \in \Int _{L_n}^{ \eta _{\la m, m^n\ra \not\hookrightarrow}}(X) $, then
$\rho \in \Max (\Int _{L_n}^{ \eta _{\la m, m^n\ra \not\hookrightarrow}}(X) )$ iff
\begin{equation}\label{EQA094}
\forall \bar x \in X^n\setminus \rho \;\;\exists K \in [X]^m \;\;K^n\setminus \rho =\{ \bar x\}.
\end{equation}
\end{cla}
\dok
(a) If there exists $\bar x\in \rho $ such that $\rho \cap K^n \setminus \{ \bar x\}\neq \emptyset$, for each $K \in [X]^m$ satisfying $\bar x \in K ^n$,
then $\rho \setminus \{ \bar x\}\in \Int _L^{ \eta _{\la m, \emptyset\ra \not\hookrightarrow}}(X) $ so, $\rho $ is not minimal.

Suppose that (\ref{EQA093}) holds and that $\rho \varsupsetneq \sigma \in \Int _L^{ \eta _{\la m, \emptyset\ra \not\hookrightarrow}}(X)$. Then, by (\ref{EQA093}),
for $\bar x \in \rho \setminus \sigma$  there is $K \in [X]^m $ such that $\rho \cap K^n =\{ \bar x\}$ and, hence, $\sigma\ cap K^n =\emptyset$, which is impossible
because $\la m, \emptyset\ra \not\hookrightarrow \la X, \sigma \ra$.

(b) follows from (a) and Theorem \ref{TA029}(c).
\kdok
\noindent
Now we show that the minimal binary structures omitting the minimal structure $\la m,\emptyset\ra$ can be characterized
using maximal $\BK_m$-free graphs.
\begin{cla}\label{TA042}
If $|X|\geq m\geq 2$, then $\rho \in \Min (\Int _{L_b}^{\eta _{\la m, \emptyset\ra \not\hookrightarrow}}(X) )$ iff $\rho $ is of the form
$$\rho=\sigma_{X\setminus R}\cup\Delta_R,$$
where $R\subset X$, $|X\setminus R|\geq m-1$, and $\sigma_{X\setminus R}$ is an orientation of the graph-complement of a maximal $\BK_m$-free graph
$\la X\setminus R, \tau_{X\setminus R} \ra$.
\end{cla}
\dok
$(\Rightarrow )$
Let $\rho \in \Min (\Int _{L_b}^{\eta _{\la m, \emptyset\ra \not\hookrightarrow}}(X) )$ and let $R:= \{x\in X: \la x,x \ra \in \rho \}$.
$|X\setminus R|\leq m-2$ would imply that for each $K\in [X]^m$ we have $|K\cap R|\geq 2$ and, hence, $|\rho\cap K^2|\geq 2$, which
is impossible by (\ref{EQA093}); so, $|X\setminus R|\geq m-1$.

By (\ref{EQA093}), for $\la x,y\ra\in \rho \cap (R\times X)$ there is $K\in [X]^m$ such that $\rho \cap K^2 =\{ \la x,y\ra\}$ and, since
$\la x,x\ra\in\rho \cap K^2$, we have $x=y$. Thus $\rho \cap (R\times X)=\Delta _R$ and, similarly, $\rho \cap (X\times R)=\Delta _R$, which means that
$\rho=\Delta_R\cup \sigma_{X\setminus R}$, where $\sigma _{X\setminus R} :=\rho \cap (X\setminus R)^2$.

By (\ref{EQA093}), for $\la x,y\ra\in \sigma _{X\setminus R}$ there is $K\in [X]^m$ such that $\rho \cap K^2 =\{ \la x,y\ra\}$ and, since
$x\neq y$, we have $\la y,x\ra\not\in \sigma _{X\setminus R}$; thus $\sigma _{X\setminus R}\cap \sigma _{X\setminus R} ^{-1} = \emptyset$.
Moreover, since
$x\neq y$, we have $K\cap R =\emptyset$, that is, $K\in [X\setminus R]^m$; so, by Claim \ref{TA041}(a)
$\sigma _{X\setminus R} \in \Min (\Int _{L_n}^{\eta _{\la m, \emptyset\ra \not\hookrightarrow}}(X\setminus R) )$.
Thus $\la X\setminus R ,\sigma _{X\setminus R} \ra$ is a minimal digraph omitting $\la m,\emptyset\ra$ and, hence,
its symmetrization $\la X\setminus R, \sigma _{X\setminus R}\cup \sigma _{X\setminus R} ^{-1}\ra$ is a minimal $\E_ m$-free graph.
By Claim \ref{TA040}(d), the graph-complement $\tau _{X\setminus R}$ of  $\sigma _{X\setminus R}\cup \sigma _{X\setminus R} ^{-1}$
is a maximal $\BK _m$-free graph and $\sigma_{X\setminus R}$ is an orientation of its graph-complement.

$(\Leftarrow )$ Let $K\in [X]^m$. If $K\cap R\neq\emptyset$, then $\rho\cap K^2\neq\emptyset$, and if $K\cap R=\emptyset$, then $\rho\cap K^2\neq\emptyset$ because otherwise the graph $\tau_{X\setminus R}$ would not be $\BK_m$-free. Hence $\rho \in \Int _{L_b}^{\eta _{\la m, \emptyset\ra \not\hookrightarrow}}(X)$.
Let $\rho ' \varsubsetneq \rho$ and $\la x,y\ra \in \rho \setminus \rho '$.
If $x=y$, take $Z\in [X\setminus R]^{m-1}$, such that $\rho\cap Z^2=\emptyset$ (such $Z$ exists because $|X\setminus R|\geq m-1$ and the graph $\tau _{X\setminus R}$ is maximal $\BK _m$-free, and thus not $\BK _{m-1}$-free by Claim \ref{TA040}(b)). Then $\rho '\cap (Z\cup \{x\})^2=\emptyset$, that is $\rho '\not \in \Int _{L_b}^{\eta _{\la m, \emptyset\ra \not\hookrightarrow}}(X)$.
If  $x\neq y$, then $x,y\in X\setminus R$, $\{x,y\} \not \in \tau _{X\setminus R}$, and since the graph $\tau _{X\setminus R}$ is maximal $\BK _m$-free, there is
$Z \subset X\setminus R$ such that $x,y\in Z $ and $\la Z, (\tau _{X\setminus R}\cup \{x,y\}) \cap Z^2 \ra \cong \BK _m$. Now $\rho '\cap Z^2 = \emptyset$ that is $\rho '\not \in \Int _{L_b}^{\eta _{\la m, \emptyset\ra \not\hookrightarrow}}(X)$. Therefore  $\rho \in \Min (\Int _{L_b}^{\eta _{\la m, \emptyset\ra \not\hookrightarrow}}(X) )$.
\kdok
In particular, for $m=2$, we have that  $\rho \in \Min (\Int _{L_b}^{\eta _{\la 2, \emptyset\ra \not\hookrightarrow}}(X) )$ if and only if there is a set $R\varsubsetneq X$ and a tournament relation $\sigma _{X\setminus R} $ on the set $X\setminus R$ such that
\begin{equation}\label{EQA019}
\rho =\sigma _{X\setminus R} \cup \Delta _R
\end{equation}
($\la X,\rho\ra$ is a {\it disjoint union of a nonempty tournament and isolated reflexive points}).

Dually we have: $\rho \in \Max (\Int _{L_b}^{ \eta _{\la 2,2^2\ra \not\hookrightarrow}}(X))$  iff there is a set $R\varsubsetneq X$ and a tournament relation
$\sigma _{X\setminus R} $ on the set $X\setminus R$ such that
\begin{equation}\label{EQA028}
\rho =X^2\setminus (\sigma _{X\setminus R} \cup \Delta _R ).
\end{equation}
For example, from (\ref{EQA019}) and (\ref{EQA028}) we obtain the reversibility of  {\it tournaments} and {\it reflexivized tournaments} (for $R=\emptyset$).
In particular, {\it strict linear orders} and {\it reflexive linear orders} are reversible.
If we take $R=X\setminus\{ x\}$, for some $x\in X$, then we obtain the {\it diagonal without one point} and {\it complete graph with one reflexivized point}.
We note that {\it complete graphs with $n$ reflexivized points} are also reversible, but for $n\geq 2$ they contain a copy of $\la 2,2^2\ra$.
\paragraph{Maximal graphs without cycle subgraphs} For $n\geq 4$, and concerning infinite graphs, the full graph is the only maximal graph which do not contain
a copy of $\C _n$, while for $n=3$ we obtain the triangle-free graphs considered above. If $3\in A\subset \o \setminus 3$, then we obtain non-trivial maximal
interpretations which do not contain copies of $\C_n$, for $n\in A$, and in the extreme situation, when we take $A=\o\setminus 3$, we obtain graphs without
cycles. The maximal such graphs are the {\it trees} (connected graphs without cycles). The {\it bipartite graphs} are obtained
if we take $A=\{ 3,5, 7, \dots\}$ and the maximal ones are the {\it complete bipartite graphs}.
\paragraph{Local cardinal bounds}
Let $L=\la R_i:i\in I\ra$ be a finite language, where $\ar (R_i)=n_i$, for $i\in I$, let $M\subset \N$ and let
$k=\la k^i_m : m\in M \land i\in I\ra$ and $l=\la l^i_m : m\in M \land i\in I\ra$ be
sequences in $\o$ such that for each $m\in M$ and $i\in I$ we have
$0\leq k^i_m \leq l^i_m \leq m^{n_i}$. Then the set of $L$-sentences
$$\textstyle
\CT _{M}^{k,l}
:=\bigcup _{m\in M}\Big\{ \eta _{\la m,\s \ra\not\hookrightarrow }:
\s \in \Int _L (m) \land \exists i\in I\;   (|\s _i|< k^i_m \lor |\s _i|> l^i_m )\Big\}
$$
is a $\Pi ^0_1$ theory and for a non-empty set $X$ and  $\r \in \Int _L (X)$ we have
$$
\r \in \Int _L ^{\CT _{M}^{k,l}}(X) \Leftrightarrow \forall m\in M \;\; \forall K\in [X]^m \;\; \forall i\in I \;\; k^i_m \leq  |\r _i \cap K^{n_i}|\leq l^i_m
$$
(the size of the components of $\r$ restricted to $m$-element subsets of $X$ is bounded).
By Theorem \ref{TA029}, if $\CT$ is an $L_{\infty \o}$-theory and the poset $\Int _L ^{\CT \cup \CT _{M}^{k,l}}(X)$ is non-empty, then it
has a dense set of minimal and co-dense set of maximal elements.
\begin{ex}\rm \label{EXA021}
Graph theory does not admit two non-trivial bounds. If $\la m, \s\ra$ is a graph, then, by irreflexivity, $0\leq |\s |\leq m^2 -m$.
Let $L=L_b$, $M=\{ 3\}$ and $0<k\leq l <6$. If
$\CT _{\{ 3\}}^{ k,l}=\{ \eta _{\la 3,\s \ra\not\hookrightarrow }: \s \in [3^{2}]^{<k}\cup [3^{2}]^{>l}\}$, then
$\CT :=\CT _{graph}\cup \CT _{\{ 3\}}^{ k,l}$
is a $\Pi ^0_1$ theory and
$\r \in \Int _{L_b} ^{\CT }(\o )$ iff the structure $\X=\la \o,\r\ra$ is a graph such that
$ k \leq  |\r  \cap K^{2}|\leq l $, for each $K\in [\o ]^3$, which means that (by symmetry)
each 3-element substructure of $\X$ has one or two edges. But this is impossible, because, by the Ramsey theorem,
$\X$ must contain an infinite empty or complete subgraph.
On the other hand, if we take $k=0$,  then for $l\in \{ 4,5\}$ the condition $ |\r  \cap K^{2}|\leq l $ means that
the graph is triangle-free and some maximal interpretations with that property are described in Examples
\ref{EXA004} and  \ref{EXA019}. For $l\in \{ 2,3\}$ maximal interpretations satisfying $ |\r  \cap K^{2}|\leq l $
are $\bigcup _\o \BK _2$ and $\BK _1 \cup \bigcup _\o \BK _2$.
\end{ex}

\section{Uniform upper bounds for definable sets}\label{S7}
If $\f (v_1,\dots ,v_p,w_1, \dots ,w_q)$ is an $L_{\infty \o }$-formula and $n\in \N$, then,
denoting $p$-tuples by $\bar {}$ and $q$-tuples by $\widetilde {}\;$, by
$\exists _{\leq n} \widetilde w \;\f (\bar v, \widetilde w)$ we denote the formula
\begin{eqnarray*}\textstyle
&&\forall w^1_1,\dots, w^1_q, \dots , w^{n+1}_1 ,\dots ,w^{n+1}_q \;\\
&&\textstyle \Big(   \bigvee _{k\leq n+1}\neg   \f (v_1,\dots ,v_p,w^k_1, \dots ,w^k_q)
\lor \bigvee _{1\leq k<l \leq n+1}\bigwedge _{1\leq j\leq q} w^k_j=w^l_j   \Big),
\end{eqnarray*}
where $\{ w^1_1,\dots, w^1_q, \dots , w^{n+1}_1 ,\dots ,w^{n+1}_q\}$ is a set of $q(n+1)$ different variables and
$\f (v_1,\dots ,v_p,w^k_1, \dots ,w^k_q)$ is the formula obtained from $\f$ by replacement of each free occurrence of $w_j$ by $w^k_j$.
Or, shortly,
$$\textstyle
\exists _{\leq n} \widetilde w \;\f (\bar v, \widetilde w):=\forall \widetilde w^1, \dots , \widetilde w^{n+1}\;
\textstyle \Big(   \bigvee _{k\leq n+1}\neg   \f (\bar v,\widetilde w^k)
\lor \bigvee _{1\leq k<l \leq n+1} \widetilde w^k=\widetilde w^l   \Big).
$$
Clearly, $\p :=\forall \bar v\; \exists _{\leq n} \widetilde w \;\f (\bar v, \widetilde w) $ is an $L_{\infty \o }$-sentence and,
if $\X$ is an $L$-structure, then $\X\models \p $ iff for each $\bar x\in X^p$ the set
$D_{\X,\f ,\bar x, q }:=\{ \widetilde y \in X^q : \X \models \f [\bar x, \widetilde y]\}$ is of size $\leq n$.
($D_{\X,\f ,\bar x, q }$ is the $q$-ary relation on the set $X$ definable in the structure $\X$ by the formula $\f$ with the parameters $x_1, \dots ,x_p$.)

Let $\neg \CF $ be the class of $L_{\infty  \o }$-formulas $\bigcup _{\xi \in \Ord}\neg \CF _\xi$, where
\begin{eqnarray*}
\neg \CF _0        &\! =  \!  &  \CN   \cup \{ R_i (v_{\a _1}, \dots , v_{\a _{n_i}}) : i\in I \land  \la \a _1, \dots , \a _{n_i} \ra \in \k ^{n_i} \},\\
\neg \CF _{\xi +1} &\! =  \!  & \neg \CF _\xi  \; \cup \; \{ \exists v_\a \; \f : \a\in \k \land  \f\in \neg \CF _\xi  \}   \\
              &\! \cup \!& \textstyle \{ \bigvee \Phi :  \Phi \subset \neg \CF _\xi  \} \;
              \cup  \; \{ \bigwedge \Phi :  \Phi \subset \neg\CF _\xi \land |\Phi |<\o \} , \\
\neg \CF _{\gamma} &\! =   \! & \textstyle  \bigcup _{\xi <\gamma }\neg\CF _\xi , \mbox{ for a limit ordinal } \gamma .
\end{eqnarray*}
Also we define the class of $L_{\infty  \o }$-formulas $\neg \CG =\bigcup _{\xi \in \Ord}\neg \CG _\xi$, where
\begin{eqnarray*}
\neg \CG _0        &\! =  \!  &  \CP  \cup \{ \neg R_i (v_{\a _1}, \dots , v_{\a _{n_i}}) : i\in I \land  \la \a _1, \dots , \a _{n_i} \ra \in \k ^{n_i} \},\\
\neg \CG _{\xi +1} &\! =  \!  & \neg \CG _\xi  \; \cup \; \{ \exists v_\a \; \f : \a\in \k \land  \f\in \neg \CG _\xi  \}   \\
              &\! \cup \!& \textstyle \{ \bigvee \Phi :  \Phi \subset \neg \CG _\xi  \} \;
              \cup  \; \{ \bigwedge \Phi :  \Phi \subset \neg\CG _\xi \land |\Phi |<\o \} , \\
\neg \CG _{\gamma} &\! =   \! & \textstyle  \bigcup _{\xi <\gamma }\neg\CG _\xi , \mbox{ for a limit ordinal } \gamma .
\end{eqnarray*}
\begin{te}\label{TA035}
Let $\f (\bar v, \widetilde w)$ be an $L_{\infty  \o }$-formula and $\CT $ an $L_{\infty  \o }$-theory such that the poset
$\P:=\la \Int _L ^{\CT \cup \{ \forall \bar v\; \exists _{\leq n} \widetilde w \;\f (\bar v, \widetilde w)\} }(X), \subset \ra$ is non-empty. Then

(a) If $\f (\bar v, \widetilde w)\in \neg \CF$ and $\CT \subset \CF$, then
$\P$ is a union-complete poset and $\Max \P $ is a co-dense set in $\P$ consisting of reversible interpretations.

(b) If $\f (\bar v, \widetilde w)\in \neg \CG$ and $\CT \subset \CG$, then
$\P$ is an intersection-complete poset and $\Min \P $ is a dense set in $\P$ consisting of reversible interpretations.
\end{te}
A proof of the theorem is given after some preliminary work.
First, to each $L_{\infty  \o }$-formula $\f$ we adjoin a formula $\f  ^\neg $ as follows:

$(v_\a =v_\b ) ^\neg  := \neg v_\a =v_\b $ and $(R_i (v_{\a _1}, \dots , v_{\a _{n_i}})) ^\neg  := \neg  R_i (v_{\a _1}, \dots , v_{\a _{n_i}}) $;

\noindent
If $\xi\in \Ord$ and $\f  ^\neg $ is defined for a formula $\f\in \Form _\xi$, then

$(\neg \f  ) ^\neg  := \f    $, $(\forall v_\a \; \f ) ^\neg  := \exists v_\a \; \f  ^\neg   $ and
$(\exists v_\a \; \f ) ^\neg  := \forall v_\a \; \f  ^\neg   $;

\noindent
If $\Phi \subset \Form _\xi$ and $\f  ^\neg $ is defined for each formula $\f\in \Phi$, then

$(\bigwedge \Phi  ) ^\neg  :=\bigvee \Phi  ^\neg    $ and $(\bigvee \Phi  ) ^\neg  := \bigwedge \Phi  ^\neg  $,

\noindent
where, for a set $\Phi$ of $L_{\infty  \o }$-formulas by $\Phi  ^\neg  $ we denote the set $\{ \f  ^\neg  : \f \in \Phi \}$.
\begin{fac}\label{TA032}
Let $\f$ be an $L_{\infty  \o }$-formula. Then

(a) $\f  ^\neg  \leftrightarrow \neg \f$;

(b) If $\f \in  \CN  $, then $\f  ^\neg \in  \CP $;

(c) If $\f \in \neg \CF $, then $\f  ^\neg \in \CF $.
\end{fac}
\dok
(a)  Let $\X\in \Mod _L$. By induction we show that for each $\f\in \Form _{L_{\infty \o }}$ we have
\begin{equation}\label{EQA033'}
\forall \vec x\in {}^{\k }X \;\Big( \X \models  \f  ^\neg  [\vec x] \Leftrightarrow  \X \models (\neg \f )[\vec x]   \Big).
\end{equation}
For $\f\in \At _{L}$ by definition we have $\f ^\neg = \neg \f$, so (\ref{EQA033'}) is true.

Let  $\f\in \Form _\xi $ and suppose that (\ref{EQA033'}) is true.
For $\vec x\in {}^{\k }X$ we have
$\X \models(\forall v_\a \; \f ) ^\neg  [\vec x] $
iff $\X\models (\exists v_\a \; \f  ^\neg  )[\vec x] $
iff for some $x\in X$ we have $\X \models \f  ^\neg  [\vec x _{\la \a , x\ra}] $,
that is, by (\ref{EQA033'}), $\X \models (\neg \f ) [\vec x _{\la \a , x\ra}] $,
that is, $\X \models \f  [\vec x _{\la \a , x\ra}] $ is not true,
iff it is not true that for all $x\in X$ we have $\X \models \f  [\vec x _{\la \a , x\ra}] $
iff $ \X \models  (\forall v_\a \; \f )  [\vec x]$ is not true
iff $ \X \models  (\neg \forall v_\a \; \f )  [\vec x]$.
Also
$\X \models(\exists v_\a \; \f ) ^\neg  [\vec x] $
iff $\X\models (\forall v_\a \; \f  ^\neg  )[\vec x] $
iff for each $x\in X$ we have $\X \models \f  ^\neg  [\vec x _{\la \a , x\ra}] $,
that is, by (\ref{EQA033'}), $\X \models (\neg \f ) [\vec x _{\la \a , x\ra}] $,
that is, $\X \models \f  [\vec x _{\la \a , x\ra}] $ is not true,
iff it is not true that for some $x\in X$ we have $\X \models \f  [\vec x _{\la \a , x\ra}] $
iff $ \X \models  (\exists v_\a \; \f )  [\vec x]$ is not true
iff $ \X \models  (\neg \exists v_\a \; \f )  [\vec x]$.
Finally,
$\X \models(\neg  \f ) ^\neg  [\vec x] $
iff $\X \models  \f [\vec x] $
iff $\X \models  \neg (\neg \f )[\vec x] $.

Let $\Phi \subset \Form _\xi $ and suppose that $\f$ satisfies (\ref{EQA033'}), for all $\f\in \Phi$.
For $\vec x\in {}^{\k }X$ we have $\X \models(\bigwedge \Phi  ) ^\neg  [\vec x] $
iff $ \X \models \bigvee \{\f  ^\neg  :\f \in  \Phi \} [\vec x] $,
iff for some $\f \in  \Phi$ we have $ \X \models \f  ^\neg  [\vec x] $,
that is, by (\ref{EQA033'}),  $ \X \models (\neg \f ) [\vec x] $,
iff it is not true that for all $\f \in  \Phi$ we have $\X \models  \f  [\vec x] $,
iff it is not true that $\X \models(\bigwedge \Phi  ) [\vec x] $
iff $\X \models(\neg \bigwedge \Phi  ) [\vec x] $.
Also,
$\X \models(\bigvee \Phi  ) ^\neg  [\vec x] $
iff $ \X \models \bigwedge \{\f  ^\neg  :\f \in  \Phi \} [\vec x] $,
iff for each $\f \in  \Phi$ we have $ \X \models \f  ^\neg  [\vec x] $,
that is, by (\ref{EQA033'}),  $ \X \models (\neg \f ) [\vec x] $,
iff it is not true that for some $\f \in  \Phi$ we have $\X \models  \f  [\vec x] $,
iff it is not true that $\X \models(\bigvee \Phi  ) [\vec x] $
iff $\X \models(\neg \bigvee \Phi  ) [\vec x] $.

(b) First we have
$(v_\a =v_\b ) ^\neg  := \neg v_\a =v_\b \in  \CP $ and $(\neg R_i (v_{\a _1}, \dots , v_{\a _{n_i}})) ^\neg  := R_i (v_{\a _1}, \dots , v_{\a _{n_i}}) \in  \CP $ and,
also, $(\neg v_\a =v_\b ) ^\neg  := v_\a =v_\b \in  \CP $.
If $\f\in \CN _\xi $ and $\f  ^\neg \in  \CP $, then the formulas
$(\forall v_\a \; \f ) ^\neg  := \exists v_\a \; \f  ^\neg   $ and
$(\exists v_\a \; \f ) ^\neg  := \forall v_\a \; \f  ^\neg   $ belong to $ \CP $.
If $\Phi \subset \CN _\xi $ and $\f  ^\neg \in  \CP $, for all $\f\in \Phi$,
then $\Phi  ^\neg \subset  \CP $ and, hence, $(\bigwedge \Phi  ) ^\neg  :=\bigvee \Phi  ^\neg   \in  \CP  $
and $(\bigvee \Phi  ) ^\neg  := \bigwedge \Phi  ^\neg  \in  \CP $.

(c) By (b) the claim is true for all $\f\in \CN  $. Also we have
$(R_i (v_{\a _1}, \dots , v_{\a _{n_i}})) ^\neg  := \neg R_i (v_{\a _1}, \dots , v_{\a _{n_i}}) \in \CF $.
If $\f\in \neg \CF _\xi $ and $\f  ^\neg \in \CF $,
then $(\exists v_\a \; \f ) ^\neg  := \forall v_\a \; \f  ^\neg   \in \CF $.
If $\Phi \subset \neg \CF _\xi $ and $\f  ^\neg \in \CF $, for all $\f\in \Phi$,
then $\Phi  ^\neg \subset \CF $ and, hence, $(\bigvee \Phi  ) ^\neg  := \bigwedge \Phi  ^\neg  \in \CF $.
If $\Phi $ is a finite set, then, again,
$(\bigwedge \Phi  ) ^\neg  :=\bigvee \Phi  ^\neg   \in \CF $.
\kdok
\noindent
{\bf Proof of Theorem \ref{TA035}.} (a) It is a standard fact that, if $\p$ is an $L_{\infty  \o }$-formula, $\f$ its subformula,
$\f \leftrightarrow \f '$ and $\p'$ the formula obtained from $\p$ by replacement of $\f$ by $\f'$,
then  $\p \leftrightarrow \p '$.
So, by Fact \ref{TA032}(a), the sentence $\p :=\forall \bar v\; \exists _{\leq n} \widetilde w \;\f (\bar v, \widetilde w) $ is
logically equivalent to the sentence
\begin{equation}\label{EQA034}\textstyle
\p ':= \forall \bar v\;\forall \widetilde w^1,\dots  \widetilde w^{n+1} \;
\Big(   \bigvee _{k\leq n+1} \f  ^\neg (\bar v, \widetilde w^k) \lor \bigvee _{1\leq k<l \leq n+1}  \widetilde w^k= \widetilde w^l   \Big).
\end{equation}
By Fact \ref{TA032}(c) we have $\f  ^\neg \in\CF$ so the sentence $\p '$ belongs to $\CF$ and the statement follows from
Theorem \ref{TA014}(a).

(b) It is easy to check that for an $L_{\infty \o }$-formula $\f$, up to logical equivalence, we have:
$\f \in  \CP $ iff $\f  ^c \in  \CN $  and, also,  $\f \in \neg \CG $ iff $\f  ^c \in \neg \CF$.
So, if $\f (\bar v, \widetilde w)\in \neg \CG$, then $\f  ^c(\bar v, \widetilde w)\in \neg \CF$ and, by (a),
$\P':=\la \Int _L ^{\CT ^c \cup \{  \forall \bar v\; \exists _{\leq n} \widetilde w \;\f   ^c (\bar v, \widetilde w)\} }(X), \subset \ra$
has the properties from (a).
Since $\forall \bar v\; \exists _{\leq n} \widetilde w \;\f   ^c (\bar v, \widetilde w)$ is the formula
$(\forall \bar v\; \exists _{\leq n} \widetilde w \;\f   ^c (\bar v, \widetilde w)) ^c$, by Theorem  \ref{TA027}(b) we have
$\P '=\{ \r ^c :\r\in \Int _L ^{\CT \cup \{ \forall \bar v\; \exists _{\leq n} \widetilde w \;\f  (\bar v, \widetilde w)\} }(X)\}$, which means that
$\P$ is the reverse of $\P '$ and, hence, has the mentioned properties.
\hfill $\Box$
\paragraph{Maximal graphs of finite degree}
If $n\in \o$, a graph $\BG=\la X, \r\ra$ is
{\it of degree $\leq n$} iff $\deg (x):=|\{ y\in X : \{ x,y \}\in \r\}|\leq n$, for all $x\in X$.
Since the atomic $L_b$-formula $R(v,w)$ belongs to the class $\neg \CF$, by Theorem \ref{TA035} the poset
$\P\!=\!\la \Int _{L_b}^{\CT _{graph}\cup \{ \p _{\deg \leq n}\}}(X), \subset \ra$, where
$\p _{\deg \leq n}:=\forall v\; \exists _{\leq n} w \; R(v,w)$ is the $L_b$-sentence
saying that a graph is of degree $\leq n$,
is union-complete and $\Max \P \subset \Rev _{L_b}(X)$ is a co-dense set in $\P$.
\begin{ex}\rm \label{EXA006}
Maximal graphs of degree $\leq 2$.  We recall that, for $n\in \o$, a graph $\BG=\la X, \r\ra$ is called
{\it $n$-regular} iff $\deg (x)=n$, for all $x\in X$.
The linear graphs $\la \o, \{ \{ n, n+1 \} : n\in \o \}\ra$ and $\la Z, \{ \{ n, n+1 \} : n\in Z \}\ra$, where $Z$ is the set of integers, will be denoted by
$\BG _\o$ and $\BG _\Z$, respectively.
\begin{cla}\label{TA017}\rm
A graph $\X$ is a maximal graph of degree $\leq 2$ iff $\X\cong \Y \;\dot{\cup }\;\Z$, where

- $\Y$ is $\emptyset$, or a 2-regular graph,

- $\Z$ is $\emptyset$, or $\BK _1$,  or $\BK _2$, or $\BG _\o$.
\end{cla}
\dok
The implication $\Leftarrow$ is evident.

($\Rightarrow$) Let $\X =\la X ,\sim \ra$ be a maximal graph of degree $\leq 2$. Suppose that
there are three different vertices $x,y,z\in X$ of degree $<2$. Then the substructure of $\X$
determined by $\{ x,y,z \}$ is not a complete graph, say $\{ x, y \}\not\in \;\sim $ and the graph
$\la X ,\sim \cup \{\{ x, y \}\}\ra $ is  of degree $\leq 2$, which contradicts the maximality of $\X$.
If all the vertices of $\X$ are of degree 2, then $\X=\Y \;\dot{\cup }\;\emptyset$, where $\Y$ is a 2-regular graph.

If two vertices of $\X$, say $x$ and $y$, are of degree $<2$, then by the maximality of $\X$, $\{ x, y \}\in \;\sim $ and, hence,
$\X=\Y \;\dot{\cup }\;\BK_2$, where $\Y$ is a 2-regular graph or $\emptyset$.

If exactly one vertex of $\X$, say $x$, is of degree $<2$ and $\deg (x)=0$, then $\X=\Y \;\dot{\cup }\;\BK_1$, where $\Y$ is a 2-regular graph or $\emptyset$.
Otherwise we have $\deg (x)=1$ and, hence, there is $y\in X\setminus \{ x \}$ such that $\{ x, y \}\in \;\sim $.
If $C_x$ is the connectivity component of $\X$ containing $x$, then in  $C_x$ we have $\deg (x)=1$ and $\deg (z)=2$, for all $z\in C_x \setminus \{ x\}$.
Now defining $x_0=x$, $x_1=y$ and $x_{n+1}$ as the unique neighbor of $x_n$ different from $x_{n-1}$ we have $\{ x_n :n\in \o\}\subset C_x$
and, by the connectedness of $C_x$ we have the equality. Thus $C_x \cong \BG _\o$.
Now, if the graph $\X$ is connected, then $\X \cong \emptyset \;\dot{\cup }\;\BG _\o$.
Otherwise, the graph induced on the set $X\setminus C_x$ is 2-regular and we have $\X\cong \Y \;\dot{\cup }\;\BG _\o$, where $\Y$ is a 2-regular graph.
\kdok
It is known that  2-regular graphs are characterized as disjoint unions of copies of $\BG _\Z$ and $\C _n$, for $n\geq 3$.
Thus there are ${\mathfrak c}$-many non-isomorphic maximal countable graphs of degree $\leq 2$; so, the poset
$\Int _{L_b}^{\CT _{graph}\cup \{ \p _{\deg \leq 2}\}}(\o )$ has ${\mathfrak c}$-many non-isomorphic maximal elements; they are reversible
and characterized in Claim \ref{TA017}.
\end{ex}
\begin{ex}\rm \label{EXA017}
Maximal connected graphs of degree $\leq n$.
Since $\f_{conn}\in  \CP $, by Theorem \ref{TA035}  maximal elements of the poset
$\la \Int _{L_b}^{\CT _{graph}\cup \{ \p _{\deg \leq n},\f_{conn}\}}(X), \subset \ra$ form a co-dense set in it consisting of reversible interpretations.
By the analysis from Example \ref{EXA006},
$\Max ( \Int _{L_b}^{\CT _{graph}\cup \{ \p _{\deg \leq 2},\f_{conn}\}}(\o ) )= [\BG _\o]_{\cong }\cup [\BG _\Z]_{\cong }$.
\end{ex}


\footnotesize
\begin{thebibliography}{99}
\bibitem{Camer}
      P.\ J.\ Cameron,
      The random graph,
      The mathematics of Paul Erd\" os II,
      Algorithms Combin. 14 (Springer, Berlin, 1997) 333--351.
\bibitem{DM}
      B. Dushnik, E.\ W.\  Miller,
      Partially ordered sets,
      Amer.\ J.\ Math.\ 63 (1941) 600--610.
\bibitem{Erdos2}
      P.\ Erd\H{o}s, A.\ R\'enyi,
      Asymmetric graphs,
      Acta Math.\ Acad.\ Sci.\ Hungar., 14 (1963) 295--315.
\bibitem{Fra1}
      R.\ Fra\"{\i}ss\'{e},
      La classe universelle: Une notion \`{a} la limite de la logique et de la th\'{e}orie des relations,
      Publ.\ D\'{e}p.\ Math.\ (Lyon) 4,3 (1967) 55--61.
\bibitem{Fra}
      R.\ Fra\"{\i}ss\'{e},
      Theory of relations, Revised edition, With an appendix by Norbert Sauer,
      Studies in Logic and the Foundations of Mathematics, 145.
      North-Holland, Amsterdam, (2000)
\bibitem{henson}
       C.\ W.\ Henson,
       A family of countable homogeneous graphs,
       Pacific J.\ Math., 38,1 (1971) 69--83.
\bibitem{Kuk}
       M.\ Kukiela,
       Reversible and bijectively related posets,
       Order 26,2 (2009) 119--124.
\bibitem{KMvar}
      M.\ S.\ Kurili\'c, N.\ Mora\v ca,
      Variations of reversibility,
      (to appear).
\bibitem{KMsim}
      M.\ S.\ Kurili\'c, N.\ Mora\v ca,
      Condensational equivalence, equimorphism, elementary equivalence and similar similarities,
      Ann.\ Pure Appl.\ Logic 168,6 (2017) 1210--1223.
\bibitem{KM3}
      M.\ S.\ Kurili\'c, N.\ Mora\v ca,
      Reversibility of disconnected structures,
      (to appear) http://arxiv.org/abs/1711.01426.
\bibitem{Lach}
       A.\ H.\ Lachlan, R.\ E.\ Woodrow,
       Countable ultrahomogeneous undirected graphs,
       Trans.\ Amer.\ Math.\ Soc., 262,1 (1980) 51--94.
\bibitem{Souk}
      L.\ Soukup,
      Infinite combinatorics: from finite to infinite,
      Horizons of combinatorics, 189--213, Bolyai Soc.\ Math.\ Stud., 17, Springer, Berlin, 2008.
\bibitem{Szp}
       E.\ Szpilrajn,
       Sur l'extension de l'ordre partiel,
       Fund.\ Math.\ 16 (1930) 386--389.
\bibitem{Tar} A.\ Tarski,
       Contributions to the theory of models. III.
       Nederl.\ Akad.\ Wetensch.\ Proc.\ Ser.\ A.\ 58 (1955) 56--64 = Indagationes Math.\ 17 (1955) 56--64.
\bibitem{Vau} R.\ L.\ Vaught,
       Remarks on universal classes of relational systems.
       Nederl.\ Akad.\ Wetensch.\ Proc.\ Ser.\ A.\ 57 = Indagationes Math.\ 16 (1954) 589--591.
\end{thebibliography}
\end{document}